\documentclass[onefignum,onetabnum]{siamart190516}

\pdfoutput=1
\usepackage{lipsum}
\usepackage{amsfonts}
\usepackage{graphicx}
\usepackage{epstopdf}
\usepackage{algorithmic}
\usepackage{xr}
\ifpdf
  \DeclareGraphicsExtensions{.eps,.pdf,.png,.jpg}
\else
  \DeclareGraphicsExtensions{.eps}
\fi
\usepackage{mwe}
\usepackage{caption}
\usepackage{subcaption}
\usepackage{listings}
\usepackage{amsmath}
\usepackage{bbm}
\usepackage{tikz}
\usetikzlibrary{positioning,shapes,angles,quotes,arrows.meta,decorations.pathreplacing,math}

\usepackage{filecontents}

\usepackage{color}
\definecolor{mygreen}{RGB}{28,172,0} 
\definecolor{mylilas}{RGB}{170,55,241}
\DeclareMathOperator{\interior}{int}
\DeclareMathOperator{\exterior}{ext}

\def\1{\mathbbm{1}}

\newsiamremark{remark}{Remark}
\newsiamremark{hypothesis}{Hypothesis}
\crefname{hypothesis}{Hypothesis}{Hypotheses}
\newsiamthm{claim}{Claim}

\headers{Fast Computation of Electrostatic Potentials}{K. Bower, K. Serkh, S. Alexakis, and A. R. Stinchcombe}

\title{Fast Computation of Electrostatic Potentials for Piecewise Constant Conductivities\thanks{Submitted to the editors 28 May 2022.
\funding{The work of the third and fourth authors was supported by the Natural Sciences and Engineering Research Council of Canada (NSERC): RGPIN-2020-01113 and RGPIN-2019-06946.}}}

\author{Kyle Bower\thanks{Department of Mathematics, University of Toronto, 40 St. George Street, Toronto ON, M5S 2E4, Canada
  (\email{kyle.bower@mail.utoronto.ca},
  \email{kserkh@math.toronto.edu}, \email{alexakis@math.toronto.edu},
  \email{stinch@math.toronto.edu}).}
  \and Kirill Serkh\footnotemark[2]
  \and Spyros Alexakis\footnotemark[2]
  \and Adam R Stinchcombe\footnotemark[2]
}

\usepackage{amsopn}

\makeatletter
\newcommand*{\addFileDependency}[1]{
  \typeout{(#1)}
  \@addtofilelist{#1}
  \IfFileExists{#1}{}{\typeout{No file #1.}}
}
\makeatother

\newcommand*{\myexternaldocument}[1]{%
    \externaldocument{#1}%
    \addFileDependency{#1.tex}%
    \addFileDependency{#1.aux}%
}

\myexternaldocument{supplement}

\ifpdf
\hypersetup{
  pdftitle={Fast Computation of Electrostatic Potentials for Piecewise Constant Conductivities},
  pdfauthor={K. Bower, K. Serkh, S. Alexakis, and A. R. Stinchcombe}
}
\fi

\begin{document}

\maketitle

\begin{abstract}
  We present a novel numerical method for solving the elliptic partial differential equation problem for the electrostatic potential with piecewise constant conductivity. We employ an integral equation approach for which we derive a system of well-conditioned integral equations by representing the solution as a sum of single layer potentials. The kernel of the resulting integral operator is smooth provided that the layers are well-separated. The fast multiple method is used to accelerate the generalized minimal residual method solution of the integral equations. For efficiency, we adapt the grid of the Nystr\"{o}m method based on the spectral resolution of the layer charge density. Additionally, we present a method for evaluating the solution that is efficient and accurate throughout the domain, circumventing the close-evaluation problem. To support the design choices of the numerical method, we derive regularity estimates with bounds explicitly in terms of the conductivities and the geometries of the boundaries between their regions. The resulting method is fast and accurate for solving for the electrostatic potential in media with piecewise constant conductivities.
\end{abstract}

\begin{keywords}
  Numerical method, Partial differential equation, Boundary integral method, Elliptic interface problem
\end{keywords}

\begin{AMS}
  65N12, 65N15, 65N22
\end{AMS}

\section{Introduction}
\label{sec:introduction}

The elliptic equation
  \begin{align}
\nabla \cdot (\sigma(\boldsymbol{x}) \nabla u(\boldsymbol{x})) = 0, \quad \boldsymbol{x}\in\Omega,
    \label{eq: PDE}
  \end{align}
in which the function $\sigma(\boldsymbol{x}) > 0$ is piecewise constant arises in numerous
important applications.  This equation has been used to model the electrical
behavior of living cells --- in which the cells and extracellular regions
possess constant conductivities~\cite{asami1980dielectric,ying_beale_2013} ---
and can be used to model the thermal
properties of composite materials~\cite{gama2017closed}.
This equation has also been used to compute the reaction potential of electrostatics when a solvation energy satisfies a nonhomogeneous jump condition for the flux~\cite{findenegg1986jn}
and to estimate the conductivity of composites formed by two isotropically conducting media~\cite{lurie1984exact}.
As the equation
modeling the electrical potential induced by a current
$\sigma(\boldsymbol{x})\partial u/\partial n$ 
injected on the boundary of some region,
it also forms
the basis of several real-world imaging techniques, known variously as
direct current resistivity (DCR), electrical capacity tomography (ECR), and
electrical impedance tomography (EIT).  The common idea behind all of these
methods is that, when a variety of currents is injected on the boundary of
some region, the electrical conductivity $\sigma(\boldsymbol{x})$ can be reconstructed
from measurements of the resulting electric potential $u(\boldsymbol{x})$, also on the
boundary (see, for example,~\cite{nachman}).

It turns out that the assumption that $\sigma(\boldsymbol{x}) > 0$ is piecewise constant
is a natural one in many applications.  In the DCR method used in
geophysics, electrode arrays are placed on the Earth's surface or inside
boreholes with the aim of determining the vertical distribution of
resistivity in the ground. The typical assumption is that geological
structures consist of pockets of various materials with differing and known
resistivities (see Table 4.3-2 in~\cite{knodel2007environmental}). 
This same imaging modality is also used to image cross-sections of
industrial processes and is known as ECT. In ECT, the capacitances of multi-electrode sensors
surrounding an industrial vessel or pipe are used to determine the makeup of
the materials inside (see, for example,~\cite{yang2002image}). Typically, the contents
consist of a multiphase flow, and the usual assumption is that there are two
materials inside the pipe or vessel with different permittivities: the flow itself, and
either inclusions of a different material or voids (bubbles)
(see~\cite{marashdeh2015electrical}).
Medical applications of this procedure go by the name of EIT, and use
sensors on the surface of the body to reconstruct an image of conductivities
inside (see, for example,~\cite{cheney}). The piecewise constant assumption is a typical one in
EIT, since it can be used to detect inclusions of differing tissues like
tumors, which have markedly different conductivities when compared with
normal tissue (see, for example,~\cite{suroviec} and~\cite{zou}).

The solution of \cref{eq: PDE} for piecewise constant $\sigma(\boldsymbol{x})$ is thus a
subject a great practical significance.  In the case of piecewise constant
conductivities, the forward problem has been solved by the finite element
method (\cite{tavares}), the finite difference method (\cite{leveque}), and
the finite volume method (\cite{dong}).  The boundary element method (BEM)
and boundary integral equation (BIEM) methods have also been applied to the
corresponding integral equation formulations (see~\cite{de-munck} and~\cite{ying_beale_2013}),
which are derived by first reformulating the forward problem as a system of
partial differential equations with transmission boundary conditions
connecting the various regions, and then representing the solution to this
system of transmission problems by layer-potentials on the boundaries of the
regions (see, for example,~\cite{colton-kress}). 
Such integral equations, when discretized, lead to well-conditioned, dense
linear systems, which can be solved rapidly using the fast multipole
method (FMM) (see, for example,~\cite{rokhlin,martinsson}).
The resulting reduction in the dimensionality of the problem and their
favorable conditioning make integral equation formulations an attractive
target for numerical methods. Boundary element methods solve these integral
equations by approximating the boundaries by polylines (see, for
example,~\cite{liu_2009}); however, when the boundaries are smooth, the number of
elements must be quite large, which reduces their overall computational
efficiency.  Boundary integral equation methods (BIEMs), on the other hand,
represent the boundaries and layer-potential densities by high order
spectral expansions, and so can require far fewer degrees of freedom, but
the evaluation of the resulting layer potentials close to their sources is
highly involved and is the subject of contemporary research
(\cite{helsing2008evaluation,barnett}).

The integral equation formulations of transmission problems have been
studied extensively from the point of view of both boundary integral
equation methods and boundary element methods, and many fast and accurate
solvers have been constructed.  However, this work has mainly focused on
acoustic and electromagnetic scattering
(\cite{gillman,hsiao,greengard-ho,perez}) and crack problems
(\cite{yoshida,otani,wang}), with some notable exceptions in the work of
Zakharov and Kalinin~\cite{zakharov} and Ying and Beale~\cite{ying_beale_2013}.  Few boundary integral
equation methods have been proposed which simultaneously address all  the
particular features of the forward electrical conductivity problem,
including: the presence of large numbers of inclusions (which may be
close-to-touching), the smooth boundaries of the inclusions, the finite
extent of the domain, the moderate accuracy requirements of the method, and
the need for reasonably accurate evaluation of the induced potential close
to and on the boundaries of the regions.

In this paper, we propose a new hybrid BEM/BIEM method for solving the
forward electrical conductivity problem (see also~\cite{seungbae}) which
constructs the solution to the equation~(\ref{eq: PDE}) on domains with
numerous inclusions, evaluates the solution quickly and accurately to within an error of
$10^{-3}$ up to the boundary, achieves $O(N)$ scaling with the use of the
FMM, and solves the problem in an automatically adaptive fashion depending
only on a user-specified accuracy.  Our hybrid method combines the efficient
layer-potential representations of BIEM with the accurate close evaluation
of BEM by interpolating between the two discretizations.  Our algorithm is
closely related to the algorithm proposed in~\cite{ying_beale_2013}, but differs in that we
show that a simpler combination of interpolation and BEM discretization 
performs just as
well as more mathematically involved methods like singularity splitting.
Since our method does not require the delicate parameter tuning common in such
methods, it is more amenable to automatic adaptivity.  We illustrate
the performance of our algorithm with several numerical examples. The method
is competitive with the state-of-the-art in the low-to-moderate accuracy
regime while being completely adaptive, and should find wide applications to
the important class of electrical conductivity problems assuming piecewise
constant conductivities.

In this paper, we also provide a completely self-contained theory, with
regards to the formulation of a well-conditioned system of boundary integral
equations, and to the existence, uniqueness, and regularity of its solutions. In
doing so, we make the following additional key contributions.  First, we
describe how the rank-deficiency resulting from the Neumann boundary
conditions of the electrical conductivity problem can be resolved by adding
an extra term to adjust the resulting system of integral equations; what is
notable is that this method does not change the dimensionality of the system
to be solved and differs from the usual approach of adjoining 
additional degrees of freedom as well as additional constraints to the system (see,
for example,~\cite{greenbaum,zakharov}). Second, we provide a detailed
justification for the choice of rescaling of the integral equations to ensure
good conditioning (see also~\cite{zakharov}). Third, we provide the proofs
of existence and uniqueness to this system of integral equations, and derive
${\cal C}^{k,\alpha}$ regularity estimates for the charge densities on the interfaces; 
these are derived using energy-type methods.
While the behavior
of such systems of boundary integral equations is considered well-understood
(see, for example,~\cite{hsiao2008boundary}), regularity results are
typically stated in the much more general setting of transmission problems
(see, for example,~\cite{costabel1988boundary}), where the dependency of the
constants on the physical parameters of the problem is often suppressed.  Our
${\cal C}^{k,\alpha}$ estimates provide bounds explicitly in terms of the conductivities
and the geometries of the boundaries between regions, and are highly
illuminating with respect to numerical and practical applications. In ``supplementary materials'', 
we also
provide stronger  estimates for the case of concentric circles, which is a
ubiquitous test case for inverse methods in EIT (see, for
example,~\cite{knudsen,cheney1990exact,mueller2003direct}). Our proofs are
also reasonably short and should be accessible to an applied audience.

This paper is organized as follows. In \cref{sec:Problem Setting}, we describe the forward electrical conductivity problem. In \cref{sec:A Boundary Integral Formulation}, we derive the corresponding system of boundary integral equations and modify it so that it is both solvable and well-conditioned. Our numerical method is presented in \cref{sec:methods}. We provide numerical experiments in~\cref{sec:experiments}. In \cref{sec:Existence and Regularity}, we prove existence and uniqueness, and provide regularity estimates for the case of general regions of constant conductivity. The conclusions follow in~\cref{sec:conclusions}.

\section{Problem Setting}
\label{sec:Problem Setting}
Let $\Omega\subset\mathbbm{R}^2$ be a simply connected and bounded domain with smooth boundary $\partial\Omega$. We consider the conductivity equation \cref{eq: PDE}
for which $u:\Omega\to\mathbbm{R}$ is the unknown scalar potential and $\sigma:\Omega\to\mathbbm{R}^+$ is the electrical conductivity. We consider the case in which the domain consists of finitely many regions of piecewise constant conductivity \textit{i}.\textit{e}., $\sigma=\sigma_i$ on region $i$ for $i=1,2,\dots,N$. The Neumann boundary condition is 
   \begin{align}
   \sigma(\boldsymbol{x})  \frac{\partial u}{\partial \boldsymbol{n}}(\boldsymbol{x}) =b_1(\boldsymbol{x}), \quad \boldsymbol{x}\in\partial\Omega,
   \label{eq: BC}
   \end{align}
in which $b_1:\partial\Omega\to\mathbbm{R}$ is the injected current. We consider the case for which $b_1\in \mathcal{C}^{0,1}$ is a Lipschitz function. We use $\boldsymbol{n}$ to denote the unit outward normal vector to the boundary and use $\frac{\partial }{\partial \boldsymbol{n}}$ or $\partial_{\boldsymbol{n}}$ for the corresponding normal derivative.

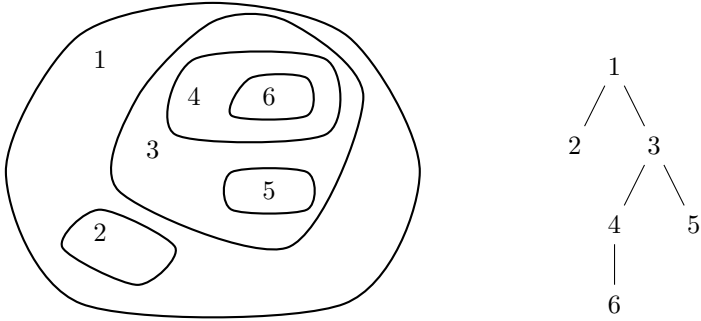
\begin{figure}[h!]
    \centering
    \begin{align*}
    \begin{tikzpicture}
    \draw[color=black] (-1.5,1.5) node {\normalsize 1};
    \draw[color=black] (-1.5,-0.8) node {\normalsize 2};
    \draw[color=black] (-0.8,0.3) node {\normalsize 3};
    \draw[color=black] (-0.25,1) node {\normalsize 4};
    \draw[color=black] (0.75,-0.25) node {\normalsize 5};
    \draw[color=black] (0.75,1) node {\normalsize 6};
    \draw [black, thick] plot [smooth cycle] coordinates {(-2.75,0) (-1.8,1.8) (0,2.25) (1.8,1.8) (2.75,0) (1.75,-1.75) (-1.75,-1.75)};
    \draw [black, thick] plot [smooth cycle] coordinates {(-1.25,-0.25) (-1,1) (0,2) (1,2) (2,1) (1,-1)};
    \draw [black, thick] plot [smooth cycle] coordinates {(-2,-1) (-1.5,-.5) (-0.5,-1) (-1,-1.5)};
    \draw [black, thick] plot [smooth cycle] coordinates {(-0.5,0.5) (-0.25,1.5) (1.5,1.5) (1.5,0.5)};
    \draw [black, thick] plot [smooth cycle] coordinates {(0.25,0.75) (0.5,1.25) (1.25,1.25) (1.25,0.75)};
    \draw [black, thick] plot [smooth cycle] coordinates {(0.25,-0.5) (0.25,0) (1.25,0) (1.25,-0.5)};
    \end{tikzpicture}
    &
    &
    \begin{tikzpicture}[every tree node/.style={draw,circle},sibling
    distance=30pt, level distance=30pt]
    \large
    \node {\normalsize 1}
        child {node {\normalsize 2}}
        child {node {\normalsize 3}
        child {node {\normalsize 4}
        child {node {\normalsize 6}}}
        child {node {\normalsize 5}}};
    \end{tikzpicture}
    \end{align*}
    \vspace{-0.6cm}
    \caption{An example layout of regions and the corresponding tree of regions. The outer boundary of the root of the tree, here region 1, is the boundary of the domain. Each region shares an interface with its children in the tree. The sets of descendants are $S_1=\{2,3,4,5,6\}, S_2=\{\}, S_3=\{4,5,6\}, S_4=\{6\}, S_5=\{\}, S_6=\{\}$.}
    \label{fig:circleTreeDiagram}
\end{figure}

We introduce a tree structure, a connected acyclic undirected graph, to describe the layout of the regions of constant conductivity. The `descendants' of a region are those regions that are entirely within it. Let $p_i$ denote the unique `parent' of region $i$ for $i=2,\ldots,N$. Let region $1$ be the root region in the tree. This is the region on which the Neumann boundary condition is imposed. The conductivity jumps across the outer boundary of region $i \ne 1$ from $\sigma_i$ to $\sigma_{p_i}$. Leaves of the tree correspond to simply connected regions that do not contain any other regions. \Cref{fig:circleTreeDiagram} shows an example layout of regions and the corresponding tree of regions. 

We seek to solve an elliptic interface problem by finding $u$ such that:
\begin{align}
    \Delta u &=0~\textrm{in}~\Omega\setminus\cup_i\partial\Omega_i \label{eq:PDE1}\\
    \sigma_1\partial_{\boldsymbol{n}} u &= b_1~\textrm{on}~\partial\Omega_1 \label{eq:PDE2}\\
    [u]&=0~\textrm{on}~\partial\Omega_i, \quad i = 2,\ldots,N \label{eq:PDE3}\\
    [\sigma\partial_{\boldsymbol{n}}u]&=b_i~\textrm{on}~\partial\Omega_i, \quad i = 2,\ldots,N \label{eq:PDE4}
\end{align}
in which $\Omega_i$ is the union of region $i$ and all its descendants. We require $\partial\Omega_i\in \mathcal{C}^2$ and $b_i\in \mathcal{C}^{0,1}$ for $i =1,\ldots,N$. For any function $\xi:\Omega\to\mathbbm{R}$ and $\boldsymbol{x}_0\in\partial\Omega_i$, we define $[\xi](\boldsymbol{x}_0)$ to be $\xi^+(\boldsymbol{x}_0)-\xi^-(\boldsymbol{x}_0)$ where $\xi^+(\boldsymbol{x}_0)$ is the limit of $\xi(\boldsymbol{x})$ as $\boldsymbol{x}\to\boldsymbol{x}_0$ from $\Omega\setminus\Omega_i$ and $\xi^-(\boldsymbol{x}_0)$ is the limit of $\xi(\boldsymbol{x})$ as $\boldsymbol{x}\to\boldsymbol{x}_0$ from $\Omega_i$. We require $\int_{\partial\Omega_i}b_i(\boldsymbol{x})~dl_{\boldsymbol{x}}=0$ for $i=1,\ldots,N$ so that the PDE problem is well-posed up to an additive constant.

The weak formulation of the PDE problem \crefrange{eq:PDE1}{eq:PDE4} will be discussed in \cref{sec:Existence and Regularity}.

\section{A Boundary Integral Formulation}
\label{sec:A Boundary Integral Formulation}

We recast the PDE problem, equations~\crefrange{eq:PDE1}{eq:PDE4}, as a boundary integral equation of the second kind~\cite{atkinson_1997}. The resulting system of integral equations is given in \cref{eq:inteqn1,eq:inteqn2}. Using an indirect approach, we represent the solution as the sum of $N$ single layer potentials
\begin{equation}
    u(\boldsymbol{x}) = \sum_{j=1}^N\mathcal{S}_{\partial \Omega_j}[\gamma_j](\boldsymbol{x}).
    \label{eq:uslpsum}
\end{equation}
The single layer potential with charge density $\gamma_j \in \mathcal{C}(\partial \Omega_j)$ is defined as
\begin{equation}
\mathcal{S}_{\partial \Omega_j}[\gamma_j](\boldsymbol{x}):=\int_{\partial \Omega_j} G(\boldsymbol{x}, \boldsymbol{y}) \gamma_j(\boldsymbol{y})~d l_{\boldsymbol{y}},
\label{eq:slp_definition}
\end{equation}
with the two-dimensional free space Green's function
\begin{equation}
    G(\boldsymbol{x}, \boldsymbol{y}) = \frac{1}{2\pi} \log|\boldsymbol{x}- \boldsymbol{y}|.
    \label{eq:freespaceGreensfunction}
\end{equation}
Since the single layer potential is harmonic off of $\partial \Omega_j$ and continuous across $\partial \Omega_j$, the sum of single layer potentials in~\cref{eq:uslpsum} automatically satisfies~\cref{eq:PDE1} and~\cref{eq:PDE3}.

As a consequence of~\cref{eq:uslpsum}, the charge density on an interface is equal to the jump in normal derivatives of the solution across the interface \textit{i}.\textit{e}.,
\begin{equation}\label{eq:charge_densities}
    \gamma_i=[\partial_{\boldsymbol{n}}u]=\partial_{\boldsymbol{n}^+}u-\partial_{\boldsymbol{n}^-}u~\textrm{on}~\partial\Omega_i, \quad i =2,\ldots,N.
\end{equation}

The normal derivative of a single layer potential on $\partial\Omega_i$ is readily computed using the adjoint of the Neumann-Poincar\'e operator,
\begin{equation}
\mathcal{K}_{\partial\Omega_i}^*[\gamma_i](\boldsymbol{x}) := \frac{1}{2\pi}\int_{\partial\Omega_i}\frac{(\boldsymbol{x}-\boldsymbol{y})\cdot\boldsymbol{n}(\boldsymbol{x})}{|\boldsymbol{x}-\boldsymbol{y}|^2}\gamma_i(\boldsymbol{y})~dl_{\boldsymbol{y}}, \quad\boldsymbol{x}\in\partial\Omega_i.
\end{equation}
A standard result from potential theory~\cite{kress2014linearintegralequations,constanda2016boundary}
is that for $\boldsymbol{x}\in\partial\Omega_i$
\begin{equation}
 \partial_{\boldsymbol{n}^\pm} \mathcal{S}_{\partial\Omega_i}[\gamma_i](\boldsymbol{x}) = \mathcal{K}_{\partial\Omega_i}^*[\gamma_i](\boldsymbol{x}) \pm \frac12\gamma_i(\boldsymbol{x}),
\end{equation}
provided that $\partial\Omega_i$ is $\mathcal{C}^2$ and $\gamma_i \in \mathcal{C}(\partial\Omega_i)$. The kernel of the adjoint of the Neumann-Poincar\'e operator,
\begin{equation}
    K(\boldsymbol{x},\boldsymbol{y})=\frac{1}{2\pi}\frac{(\boldsymbol{x}-\boldsymbol{y})\cdot\boldsymbol{n}(\boldsymbol{x})}{|\boldsymbol{x}-\boldsymbol{y}|^2},\label{eq:K_ij_kernel}
\end{equation}
has a removable discontinuity at $\boldsymbol{x}=\boldsymbol{y} \in \partial\Omega_i$. This discontinuity is removed by defining
\begin{equation}
    K(\boldsymbol{x},\boldsymbol{x}):=\frac{\kappa(\boldsymbol{x})}{4\pi},\label{eq:K_ij_kernel_removable_discontinuity}
\end{equation}
in which $\kappa$ is the curvature at point $\boldsymbol{x} \in \partial\Omega_i$. Since the interfaces are in $\mathcal{C}^2$, do not self-intersect, and do not intersect with other interfaces, $\boldsymbol{x}=\boldsymbol{y}$ implies $\boldsymbol{n}(\boldsymbol{x})=\boldsymbol{n}(\boldsymbol{y})$ and $K(\boldsymbol{x},\boldsymbol{y})\in \mathcal{C}^2$.

On the boundary of each region $\partial\Omega_i$, there is a single layer of charge with density $\gamma_i$. Enforcing the boundary condition~\cref{eq:PDE2} and the interfaces conditions~\cref{eq:PDE4} provide a system of integral equations that determine $\gamma_1,\dots,\gamma_N$. The normal derivative of the single layer potential $\mathcal{S}_{\partial\Omega_j}[\gamma_j]$ at point $\boldsymbol{x} \in \partial\Omega_i, i\in\{1,\dots,N\}$ is
\begin{equation}
\partial_{\boldsymbol{n}^\pm}\mathcal{S}_{\partial\Omega_j}[\gamma_j](\boldsymbol{x}) = \int_{\partial\Omega_j} K(\boldsymbol{x},\boldsymbol{y}) \gamma_j(\boldsymbol{y})~dl_{\boldsymbol{y}} \pm \frac12 \gamma_j (\boldsymbol{x})\delta_{ij},
\label{eq:slpnormderiv}
\end{equation}
in which $\delta_{ij}$ is the Kronecker delta.
Substituting~\cref{eq:uslpsum} and~\cref{eq:slpnormderiv} into~\cref{eq:PDE2} and~\cref{eq:PDE4} gives a system of integral equations for the charge densities,
\begin{align}
 -\frac12 \sigma_1 \gamma_1(\boldsymbol{x}) + \sigma_1 \sum_{j=1}^N \int_{\partial\Omega_j} K(\boldsymbol{x},\boldsymbol{y}) \gamma_j(\boldsymbol{y})~dl_{\boldsymbol{y}} &= b_1(\boldsymbol{x}),\quad \boldsymbol{x}\in\partial\Omega_1,\label{eq:integral_equation_1}\\
 \frac12\left(\sigma_{p_i}+\sigma_i\right)\gamma_i(\boldsymbol{x}) + \left(\sigma_{p_i}-\sigma_i\right)\sum_{j=1}^N\int_{\partial\Omega_j} K(\boldsymbol{x},\boldsymbol{y}) \gamma_j(\boldsymbol{y})~dl_{\boldsymbol{y}} &= b_i(\boldsymbol{x}), \quad \boldsymbol{x}\in\partial\Omega_i, \label{eq:integral_equation_2} \\
 &\mathrel{\phantom{=}}\text{for}~i =2,\ldots,N.\nonumber
\end{align}
The system of equations \cref{eq:integral_equation_1,eq:integral_equation_2} is ill-posed, owing exclusively to the well-known non-uniqueness for the Neumann problem for equations of this type. In particular, solutions are only defined up to an additive constant. 

 We will resolve this issue here, by making a suitable change to the system. Prior to doing this, it is instructive to derive that the total charge on any inner interface must vanish as a consequence of the system \cref{eq:integral_equation_1,eq:integral_equation_2}. The non-uniqueness is entirely captured in the freedom to define a weighted integral of the charge density on the outer boundary arbitrarily, as we will see later.

Denote the total charge on $\partial\Omega_i$ by $C_i$ \textit{i}.\textit{e}.,
\begin{equation}
    C_i:=\int_{\partial\Omega_i}\gamma_i(\boldsymbol{x})~dl_{\boldsymbol{x}},\quad\text{for}~ i =1,\ldots,N.
\end{equation}
We will show that the integral equations \cref{eq:integral_equation_2} imply that $C_i=0$ for all $i =2,\ldots,N$. We integrate the integral equations \cref{eq:integral_equation_2} along $\partial\Omega_i$ and use the requirement that $\int_{\partial\Omega_i}b_i(\boldsymbol{x})~dl_{\boldsymbol{x}} =0$ for $i\in\{2,\dots,N\}$ to obtain the system of equations
\begin{align}\label{eq:integrate_integral_equation2}
    0 &=\frac{1}{2}(\sigma_{p_i}+\sigma_i)C_i+ \left(\sigma_{p_i}-\sigma_i\right)\sum_{j=1}^N\int_{\partial\Omega_i}\int_{\partial\Omega_j} K(\boldsymbol{x},\boldsymbol{y}) \gamma_j(\boldsymbol{y})~dl_{\boldsymbol{y}}~dl_{\boldsymbol{x}}, \\
    &\mathrel{\phantom{=}}\text{for}~ i =2,\ldots,N.\nonumber
\end{align}
Since $\boldsymbol{y} \in \partial\Omega_j$, the divergence theorem and in the case when $\boldsymbol{y}\in\partial\Omega_i$ slightly modifying the contour can be used to compute the integral
\begin{align}\label{eq:integral_computed_using_divergence_theorem}
    \int_{\partial\Omega_i}K(\boldsymbol{x},\boldsymbol{y})~dl_{\boldsymbol{x}} =
    \begin{cases}
      1, & \boldsymbol{y}\in\interior(\Omega_i) \\
      \frac{1}{2}, & \boldsymbol{y}\in\partial\Omega_i \\
      0, & \boldsymbol{y}\in\exterior(\Omega_i)
   \end{cases}.
\end{align}

Swapping the order of integration in \cref{eq:integrate_integral_equation2} and using \cref{eq:integral_computed_using_divergence_theorem} leads to the system of equations
\begin{equation}
    0 = \sigma_{p_i}C_i+(\sigma_{p_i}-\sigma_i)\sum_{j\in S_{i}}C_j,
    \quad\text{for}~ i =2,\ldots,N,
\end{equation}
in which $S_{i}\subseteq\{1,\dots,N\}$ is the set of all regions that are descendants of region $i$ as illustrated in \cref{fig:circleTreeDiagram}. These equations immediately give us that for any region $i$ that has no descendants \textit{i}.\textit{e}., a leaf of our tree of regions, $C_i=0$. Then, for any parent region $p_i$ whose descendants all have $C_i=0$, we must also have $C_{p_i}=0$. By induction, we get $C_i=0$ for all $i = 2,\ldots,N$.

Subsequently, we will show that the freedom in selecting the total charge on the outer boundary is equivalent to choosing the additive constant for the PDE problem. We will begin by showing that the freedom in selecting a weighted integral of the charge density on the outer boundary is equivalent to choosing the additive constant for the PDE problem. Assume we have two solutions to the system of integral equations \cref{eq:integral_equation_1,eq:integral_equation_2}. We will denote these two solutions by $u_1$ and $u_2$. Define $u$ to be the difference between the two solutions \textit{i}.\textit{e}., $u=u_1-u_2$. 

Then $u$ satisfies the same problem as before~\crefrange{eq:PDE1}{eq:PDE4} with all the right hand sides being zero. Any solution to the system~\crefrange{eq:PDE1}{eq:PDE4} is a solution of the weak formulation of the elliptic interface problem. The only solution to the weak problem with all the right hand sides being zero is the constant solution $u=d$. 

Let us denote the charge densities corresponding to the difference of the two solutions to the integral equations system \cref{eq:integral_equation_1,eq:integral_equation_2} by  $\delta\gamma_i, i=1,\dots, N$. Since we know that the solution $u$ must be a constant, 
we have using \cref{eq:charge_densities} that $\delta\gamma_i=0$ for all $i\ge 2$. 
Therefore, using \cref{eq:uslpsum,eq:slp_definition} $u$ satisfies the formula 
\begin{equation}
    u(\boldsymbol{x})=\int_{\partial\Omega_1}G(\boldsymbol{x}, \boldsymbol{y}) \delta\gamma_1(\boldsymbol{y})~d l_{\boldsymbol{y}},
\end{equation}
in which $\delta\gamma_1$ is the charge density on the outer boundary. Integrating $u$ along $\partial\Omega_1$ and swapping the order of integration gives
\begin{equation}
    d\cdot\textrm{length}(\partial\Omega_1) =\int_{\partial\Omega_1} w(\boldsymbol{y}) \delta\gamma_1(\boldsymbol{y})~d l_{\boldsymbol{y}},
\end{equation}
in which $w(\boldsymbol{y}):=\int_{\partial\Omega_1}G(\boldsymbol{x}, \boldsymbol{y})~d l_{\boldsymbol{x}}$. If $\partial\Omega_1$ is a circle of radius $R$, then $w(\boldsymbol{y})$ is constant and identically equal to $R\log{R}$.

Next we will show that specifying the given weighted integral of the charge density on the outer boundary is equivalent to specifying the total charge on the outer boundary. We note that $\delta\gamma_1$ is the charge density on the boundary of a conductor and is thus continuous and never zero~\cite{petrovsky2012lectures}. Hence, $\delta\gamma_1$ is of constant sign. Now we have already shown there is a one-dimensional kernel to our system of integral equations \cref{eq:integral_equation_1,eq:integral_equation_2} consisting of charge densities that yield the constant solution so specifying the total charge on the outer boundary will uniquely specify the weighted integral of the charge density on the outer boundary and vice versa.

We thus derive that any two solutions to the original problem \crefrange{eq:PDE1}{eq:PDE4} with the same $b_i$ for $i=1,\dots,N$ differ by a constant and that specifying the total charge on the outer boundary uniquely specifies a solution to problem \crefrange{eq:PDE1}{eq:PDE4}. 

For definiteness, we will be imposing that the total charge on the outer boundary $C_1:=\int_{\partial\Omega_1} \gamma_1(\boldsymbol{y})~d l_{\boldsymbol{y}}=0$. We incorporate our choice of $C_1=0$ into the system of linear integral equations given in \cref{eq:integral_equation_1,eq:integral_equation_2} by modifying the kernel in \cref{eq:K_ij_kernel} when $i=j=1$ by subtracting $1$. The integral equations \cref{eq:integral_equation_1,eq:integral_equation_2} for the charge densities with the modified kernel where $C_1=0$ become
\begin{align}
 -\frac12 \sigma_1 \gamma_1(\boldsymbol{x}) + \sigma_1 \sum_{j=1}^N \int_{\partial\Omega_j}\left[K(\boldsymbol{x},\boldsymbol{y})-\delta_{1j}\right]\gamma_j(\boldsymbol{y})~dl_{\boldsymbol{y}} &= b_1(\boldsymbol{x}),\quad \boldsymbol{x}\in\partial\Omega_1,\label{eq:integral_equation_1_modified_K}\\
 \frac12\left(\sigma_{p_i}+\sigma_i\right)\gamma_i(\boldsymbol{x}) + \left(\sigma_{p_i}-\sigma_i\right)\sum_{j=1}^N\int_{\partial\Omega_j}K(\boldsymbol{x},\boldsymbol{y})\gamma_j(\boldsymbol{y})~dl_{\boldsymbol{y}} &= b_i(\boldsymbol{x}), \quad \boldsymbol{x}\in\partial\Omega_i, \nonumber \\
 &\mathrel{\phantom{=}}\text{for}~i = 2,\ldots,N\label{eq:integral_equation_2_modified_K} .
\end{align}

We seek to show that the previous integral equations  \cref{eq:integral_equation_1,eq:integral_equation_2} with the additional constraint $C_1:=\int_{\partial\Omega_1}\gamma_1(\boldsymbol{y})~dl_{\boldsymbol{y}}=0$ are equivalent to the modified integral equations \cref{eq:integral_equation_1_modified_K,eq:integral_equation_2_modified_K}. We first show that a solution to the previous integral equations  \cref{eq:integral_equation_1,eq:integral_equation_2} where in addition $C_1=0$ yields a solution to the modified integral equations \cref{eq:integral_equation_1_modified_K,eq:integral_equation_2_modified_K}. Since the only difference between the two systems of integral equations is that \cref{eq:integral_equation_1_modified_K} has an additional term on the left hand side
\begin{equation}\label{eq:integral_equation_1_modified_K_additional_term}
    -\int_{\partial\Omega_1}\gamma_1(\boldsymbol{y})~dl_{\boldsymbol{y}}
\end{equation}
which is 0 when $C_1=0$, a solution to the previous integral equations  \cref{eq:integral_equation_1,eq:integral_equation_2} where in addition $C_1=0$ yields a solution to the modified integral equations \cref{eq:integral_equation_1_modified_K,eq:integral_equation_2_modified_K}.

Next, we show the converse that a solution to the modified integral equations \cref{eq:integral_equation_1_modified_K,eq:integral_equation_2_modified_K} yields a solution to the previous integral equations  \cref{eq:integral_equation_1,eq:integral_equation_2} where in addition $C_1=0$. Integrating the integral equation with the modified kernel \cref{eq:integral_equation_1_modified_K} around $\partial \Omega_1$, swapping the order of integration, and using the divergence theorem immediately gives $    -\int_{\partial\Omega_1}\gamma_1(\boldsymbol{y})~dl_{\boldsymbol{y}}=0$
which by definition of $C_1$ implies that $C_1=0$. Now, the additional term \cref{eq:integral_equation_1_modified_K_additional_term} in equation \cref{eq:integral_equation_1_modified_K} vanishes so a solution to the modified integral equations \cref{eq:integral_equation_1_modified_K,eq:integral_equation_2_modified_K} yields a solution to the previous integral equations  \cref{eq:integral_equation_1,eq:integral_equation_2}.

Thus, the system of modified integral equations \cref{eq:integral_equation_1_modified_K}-\cref{eq:integral_equation_2_modified_K} is equivalent to the original system of integral equations \cref{eq:integral_equation_1}-\cref{eq:integral_equation_2} with the additional constraint that the total charge on the outer boundary $C_1=0$.

\subsection{Conditioning of Linear System}
\label{sec:Conditioning of Linear System}

The scenarios we are interested in typically involve conductivities of order 1; however, if two neighboring conductivities are very similar in value then the resulting system of integral equations will be poorly conditioned. For an inner interface (\textit{i}.\textit{e}., $i\ge2$), if $b_i\equiv[\sigma\partial_{\boldsymbol{n}}u]=0$, then \cref{eq:PDE4,eq:charge_densities} imply that the charge density on interface $\partial\Omega_i$ is
\begin{equation}
    \gamma_i=\partial_{\boldsymbol{n}^+}u\left(1-\frac{\sigma_{p_i}}{\sigma_i}\right)\label{eq:gammai}.
\end{equation}
This quantity will be very small if $\sigma_{p_i}\approx\sigma_{i}$ leading to a poorly conditioned system of integral equations. To improve the conditioning of the system of integral equations, we rescale the charge densities to be proportional to $\partial_{\boldsymbol{n}^+}u$ by defining $\phi_i := \gamma_i/\alpha_i$ where 
$\alpha_i = 1-\sigma_{p_i}/\sigma_i$ 
for $i\in \{2,\ldots,N\}$ and $\alpha_1 = 1$. The integral equations for the rescaled charge densities $\phi_i$ become
\begin{align}
 -\frac12 \phi_1(\boldsymbol{x}) + \sum_{j=1}^N \alpha_j \int_{\partial\Omega_j}\left[K(\boldsymbol{x},\boldsymbol{y})-\delta_{1j}\right]\phi_j(\boldsymbol{y})~dl_{\boldsymbol{y}} &= \frac{b_1(\boldsymbol{x})}{\sigma_1},\quad \boldsymbol{x}\in\partial\Omega_1,\label{eq:inteqn1}\\
 \frac12\alpha_i\frac{\left(\sigma_{p_i}+\sigma_i\right)}{\sigma_{p_i}-\sigma_i}\phi_i(\boldsymbol{x}) + \sum_{j=1}^N\alpha_j\int_{\partial\Omega_j}K(\boldsymbol{x},\boldsymbol{y})\phi_j(\boldsymbol{y})~dl_{\boldsymbol{y}}&= \frac{b_i(\boldsymbol{x})}{\sigma_{p_i}-\sigma_i},\quad \boldsymbol{x}\in\partial\Omega_i,  \label{eq:inteqn2}\\
 &\mathrel{\phantom{=}}\text{for}~ i = 2,\ldots,N\nonumber
\end{align}
in which $\delta_{1j}$ is the Kronecker delta and $K$ is the kernel of the adjoint of the Neumann-Poincar\'e operator defined in \cref{eq:K_ij_kernel,eq:K_ij_kernel_removable_discontinuity}. Thus, $\phi_i$ in view of \cref{eq:gammai} will be proportional to $\partial_{\boldsymbol{n}^+}u$ leading to a well-conditioned system of integral equations.

\section{Methods}
\label{sec:methods}
This section describes the discretization scheme for the integral equations, specifies the interpolation of the charge densities, documents the adaptive quadrature method used to automatically refine the quadrature grids, details the solution evaluation, and presents the exact solution for the two non-concentric nested circles case.

\subsection{Discretization Scheme}

We employ the Nystr\"{o}m method~\cite{kress2014linearintegralequations} to discretize the integral equations~\cref{eq:inteqn1,eq:inteqn2}.
The curves $\partial \Omega_i$ are parameterized in $q\in[0,1]$ and the grid on each curve is denoted $q_{i,m},~m=0,\ldots,M_i-1$ for $M_i$ grid points on curve $i$. We use two different quadrature schemes. For interfaces that approach other interfaces within a threshold, we use composite quadrature and split the interface into panels. We use Gauss-Legendre quadrature on each panel and Lagrange interpolation for interpolating the charge density. For interfaces that are sufficiently far away from all other interfaces we use uniform grid points, trapezoidal rule quadrature, and trigonometric interpolation for interpolating the charge density. For either interface type, we adaptively refine the panels and increase the number of grid points to meet a given error tolerance. 

After discretization, the integral equations \cref{eq:inteqn1,eq:inteqn2} can be expressed as $A\boldsymbol{\phi} = \boldsymbol{b}$. The entries in $\boldsymbol{b}$ are given by evaluating the functions $b_1,b_2,\ldots,b_N$ at the quadrature nodes and the entries in $A$ are computed from evaluating the kernel \cref{eq:K_ij_kernel} at the pairs of grid points. 

We solve the system $A\boldsymbol{\phi} = \boldsymbol{b}$ using GMRES~\cite{saad1986gmres}, an iterative method that approximates the solution by the vector in a Krylov subspace with minimal residual. When solving with GMRES, the well-conditioned integral equations give a well-conditioned discrete linear system for which GMRES converges quickly.

We use the \texttt{fmm2d} library that implements fast multipole methods in two dimensions~\cite{fmm2d}. It allows us to compute $N$-body interactions governed by the Laplace equation, to a specified precision, in two dimensions, on a multi-core shared-memory machine. In particular, \texttt{fmm2d} can accelerate evaluation of the sum
\begin{equation}
    \sum_{j=1}^M c_j \mathcal{K}\left(\boldsymbol{x}_i,\boldsymbol{y}_j\right),\quad\textrm{for}\ i=1,2\dots,N,
\end{equation}
in which there are $M$ arbitrary source locations $\boldsymbol{y}_j\in\mathbbm{R}^2$ with corresponding strengths $c_j\in\mathbbm{R}$, $N$ target locations $\boldsymbol{x}_i\in\mathbbm{R}^2$, and the kernel $\mathcal{K}$ can be the free-space Green's function $G(\boldsymbol{x}, \boldsymbol{y})$ from \cref{eq:freespaceGreensfunction} or its derivative $K(\boldsymbol{x}, \boldsymbol{y})$ from \cref{eq:K_ij_kernel}.

\subsection{Interpolation of the Charge Densities}
The boundary integral method results in estimates of the charge densities on a grid $\gamma_i(q_{i,m})$. When evaluating the singular layer potentials 
$\mathcal{S}_{\partial\Omega_i}[\gamma_i](\boldsymbol{x})=\int_{\partial\Omega_i} G(\boldsymbol{x}, \boldsymbol{y}) \gamma_i(\boldsymbol{y}) d l_{\boldsymbol{y}}$, the value of the charge density may need to be known off of this grid depending on the quadrature scheme employed.

While we could have opted to use Nystr\"{o}m interpolation~\cite{nystrom1930praktische,delves1988computational,kress2014linearintegralequations} in which the integral equations are re-arranged and evaluated for the charge density at non-grid points, we instead choose barycentric Lagrange interpolation~\cite{berrut2004barycentric} for the paneled interfaces and trigonometric interpolation~\cite{wright2015extension} for the uniform interfaces.

Consider the case of a uniform interface first. Let $p_K(q)$ be the trigonometric polynomial of degree $K$, \textit{i}.\textit{e}., a truncated Fourier series. The $2K+1$ unknowns are determined by requiring that the interpolant pass through given values $p_K(q_m)=f_m, m=0,\ldots,M-1$ with $M=2K+1$. In the case of the shifted uniform grid $q_m = (m+s)/M,~m=0,\ldots,M-1,~s\in[0,1]$, this interpolant can readily and stably be evaluated using the barycentric formula for trigonometric interpolation~\cite{wright2015extension},
\begin{equation}
p_K(q) = \frac{\displaystyle \sum_{m=0}^{M-1} (-1)^m f_m F\left(\pi(q-q_m)\right) } {\displaystyle \sum_{m=0}^{M-1} (-1)^m F\left(\pi(q-q_m)\right) }.
\end{equation}
For odd $M$, $F(q) = \csc(q)$, while for even $M$, $F(q) = \cot(q)$. We choose $s=\frac12$ so as to avoid a possible numerical instability with the barycentric formula~\cite{austin2017numerical}. 

In the case of a non-uniform grid $q_m, m=0,\ldots,M-1$, we opt for a degree $M-1$ Lagrange interpolant using the barycentric formula~\cite{berrut2004barycentric},
\begin{equation}
p_{M-1}(q) = \frac{\displaystyle \sum_{m=0}^{M-1}\frac{w_mf_m}{q-q_m}}{\displaystyle \sum_{m=0}^{M-1}\frac{w_m}{q-q_m}},
\end{equation}
in which
\begin{equation}
w_m = \frac{1}{\prod_{n \ne m} (q_m - q_n)}.
\end{equation}
If the coefficients $w_m$ are pre-computed, then interpolation requires only $\mathcal{O}(NM)$ work for $N$ evaluation points. 

\subsection{Adaptive Grid Refinement}
We use adaptive quadrature to automatically refine the quadrature grids for the integral equations.

For interfaces discretized with uniform grids, we solve the discretized integral equations~\cref{eq:inteqn1,eq:inteqn2} with $M$ uniform grid points and use trigonometric interpolation to approximate the scaled charge densities. We take $M=2K+1$ and let $p_K(q)$ be the truncated Fourier series of degree $K$ approximation to the discrete charge densities
\begin{equation}
    p_K(q)=\sum_{k=-K}^Kc_k\exp(2\pi ikq),
\end{equation}
in which the Fourier coefficients are
\begin{equation}
    c_k=\int_0^1p_K(q)\exp(-2\pi ikq)dq.
\end{equation}
If the two highest mode Fourier coefficients $|c_K|$ and $|c_{K-1}|$ are both smaller than a specified tolerance, then the charge densities are well resolved and $p_K(q)$ is taken to be our final approximation of the scaled charge density. Otherwise, we refine the interface and iterate until $|c_K|$ and $|c_{K-1}|$ are both smaller than the specified tolerance. We check that the last two coefficients, rather than just one, are smaller than the threshold so that a symmetry does not mislead us into thinking that the series has converged.

For interfaces discretized with composite quadrature (panels), we solve the discretized integral equations~\cref{eq:inteqn1,eq:inteqn2} with an $M$ point Gauss-Legendre quadrature and examine the truncated Legendre series polynomial of order $M$ approximation to the discrete charge densities. If either of the last two coefficients in the polynomial approximation are above a specified tolerance, then the panel is refined. After refinement, if the number of quadrature nodes is greater than the maximum allowed on a panel, then we split the panel into four smaller panels. We iterate until the last two coefficients in the polynomial approximation are both below the specified tolerance, indicating that the charge densities are well resolved. Refining based on the size of the coefficients of Legendre series polynomials has been used previously when computing generalized Gaussian quadratures~\cite{bremer2010generalizedGauss}.

Each time the grid is refined, the linear system for charge densities needs to be solved again. We reduced the number of GMRES iterations needed to solve this linear system by using an initial guess from interpolating the charge densities from the course grid to the fine grid.

\subsection{Solution Evaluation}
Having solved for the charge densities, the solution may be evaluated using~\cref{eq:uslpsum} and the same quadrature scheme used in the Nystr\"{o}m method. This is known as na\"{i}ve quadrature and is fast and accurate when the evaluation point is not near or on an interface or the boundary. For evaluation points on or near an interface or the boundary, the integrand is logarithmically singular or nearly singular for which the quadrature schemes for smooth integrands perform poorly~\cite{barnett2014evaluation,helsing2008evaluation}.

This so-called close evaluation problem has been addressed with quadratures that are spectrally accurate up to the boundary~\cite{helsing2008evaluation,barnett2015spectrally}, quadrature by expansion~\cite{klockner2013qbx} (whose convergence is analyzed in~\cite{epstein2013convergence}), and adaptive quadrature~\cite{gonnet2009adaptive}. We opt for a simpler approach: approximate interfaces or the boundary with line segments near an evaluation point and employ an analytic expression for the single layer potential due to a line segment.

Through explicit integration, the potential at point $\boldsymbol{x}$ due to a single layer of charge on a line segment $\Delta S$, parameterized by arclength $s$, from $\boldsymbol{X}_1$ at $s_1$ to $\boldsymbol{X}_2$ at $s_2$ with a charge density $\gamma(s) = \gamma_1 \frac{s_2-s}{s_2-s_1} + \gamma_2 \frac{s-s_1}{s_2-s_1}$ varying linearly from $\gamma_1$ at $\boldsymbol{X}_1$ to $\gamma_2$ at $\boldsymbol{X}_2$ is
\begin{align}
    \mathcal{S}_{\Delta S}[\gamma(s)](\boldsymbol{x}) = &\frac{1}{2\pi}\frac{\gamma_1 s_2 - \gamma_2 s_1}{s_2-s_1}\left(s_2\log(r_2) - s_1\log(r_1) - s_2 + s_1 +(\theta_2-\theta_1)d\right)  \label{eq:linesegmentchargepotential}\\
    & + \frac{1}{2\pi}\frac{\gamma_2-\gamma_1}{s_2-s_1}\left(\frac{r_2^2}{2}\left(\log(r_2)-\frac12\right) - \frac{r_1^2}{2}\left(\log(r_1)-\frac12\right)\right)\nonumber,
\end{align}
in which: $r_1,r_2$ are the distances between $\boldsymbol{x}$ and $\boldsymbol{X}_1,\boldsymbol{X}_2$, respectively; $d$ is the perpendicular distance between $\boldsymbol{x}$ and the line segment $\Delta S$; and $\theta_1,\theta_2$ are the angles from the perpendicular line to $\boldsymbol{X}_1,\boldsymbol{X}_2$, respectively. The variables in this formula are shown in \cref{fig:linechargepotential_a}. This formula is given in the Appendix A of~\cite{liu_2009} in the case of $\gamma_1=\gamma_2=1$.

The values of the charge densities at the endpoints of each line segment can be obtained by interpolating the charge density solution. If we simply take the values of the unmodified charge densities at the endpoints of the curved segment, then the total charge for the line segment approximation would be reduced as the length of the line segment is shorter than that of the curved segment. To ensure the total charge on each line segment is equal to the total charge on the corresponding curved segment, we estimate the total charge on each curved segment by computing a truncated Legendre (respectively Fourier) series polynomial for panelled (respectively uniform) interfaces. Imposing that the  estimated total charge on each curved segment equals that on the corresponding line segment yields a system for the charge densities at each endpoint. For an even number of endpoints, this system has a one-dimensional kernel and we choose the solution with minimum $\ell_2$-norm.

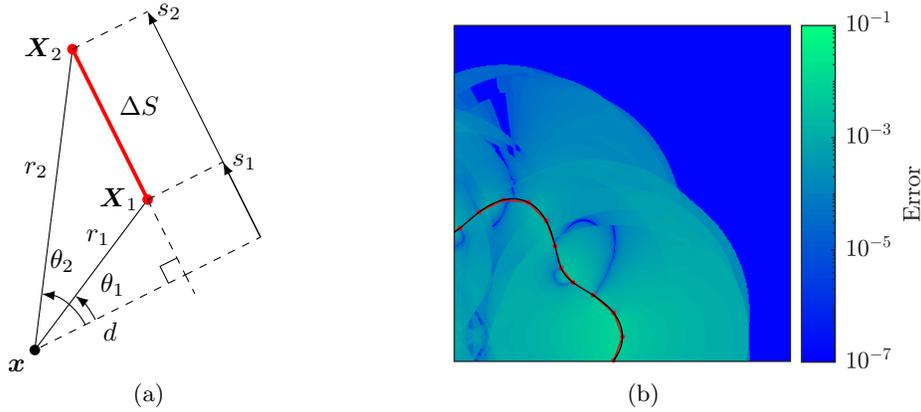
\begin{figure}
    \begin{subfigure}[b]{0.40\textwidth}
    \hspace{5mm}
    \begin{tikzpicture}
        \draw (0,0) node[left] {$\boldsymbol{X}_1$} -- node[midway, above right] {$\Delta S$} (-1,2) node[left] {$\boldsymbol{X}_2$} -- node[pos=0.4, left] {$r_2$} (-1.5,-2) node[below left] {$\boldsymbol{x}$} -- node[pos=0.75, left] {$r_1$} cycle; 
        \draw[line width =0.5mm,red] (0,0) -- (-1,2); 
        \fill[red] (0,0)  circle[radius=2pt];
        \fill[red] (-1,2)  circle[radius=2pt];
        \fill (-1.5,-2)  circle[radius=2pt];
        \draw[dashed] (0,0) -- node[right] {} (0.5,-1) -- node[midway, below] {$d$} (-1.5,-2); 
        \draw[dashed] (0.5,-1) -- (1.5,-0.5);
        \draw[dashed] (0,0) -- (1,0.5);
        \draw[dashed] (-1,2) -- (0,2.5);
        \draw[dashed] (0.5,-1) -- (0.625,-1.25);
        \draw(0.3882,-0.7764)--(0.1646,-0.8882)-- (0.2764,-1.1118);
        \draw[-Latex] (1.5,-0.5) -- (1,0.5) node[right] {$s_1$};
        \draw[-Latex] (1.5,-0.5) -- (0,2.5) node[right] {$s_2$};
         \draw[opacity=0] (0.5,-1) coordinate (P) node {} -- node[midway, below] {$d$} (-1.5,-2) coordinate (x) node {} -- (0,0) coordinate (X1) node {} -- (-1,2) coordinate (X2) node {}
         pic["$\theta_1$", draw=black, opacity=1, -Latex, angle eccentricity=1.5, angle radius=0.90cm]{angle=P--x--X1}
         pic["$\theta_2$",above left, draw=black, opacity=1, -Latex, angle eccentricity=1.5, angle radius=0.75cm]{angle=P--x--X2};
    \end{tikzpicture}
        \vspace{-0.2cm}
    \caption{}
    \label{fig:linechargepotential_a}
    \end{subfigure}
    \begin{subfigure}[b]{0.59\textwidth}
    \hspace{10mm}
    \includegraphics[]{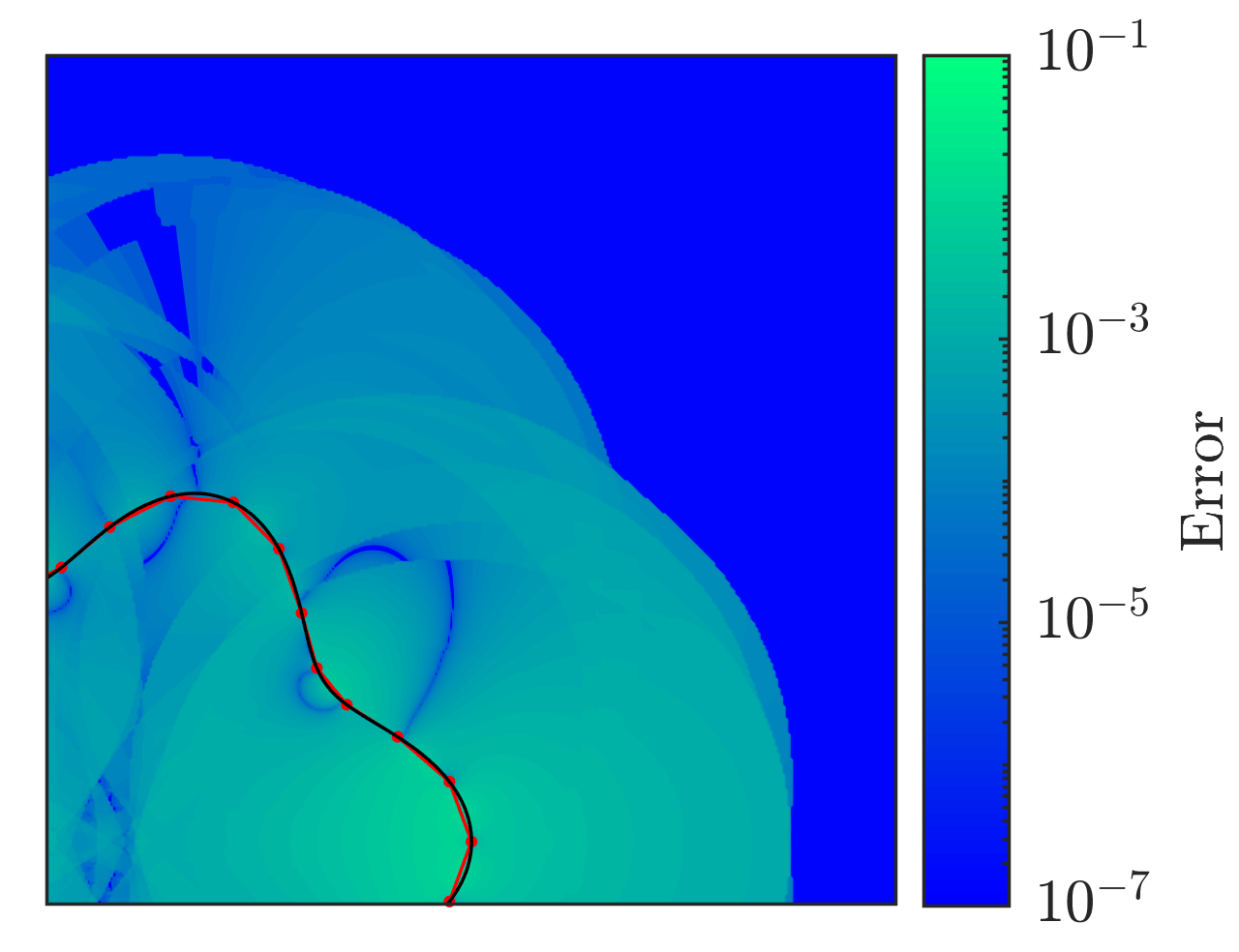}
        \vspace{-0.2cm}
    \caption{}
    \label{fig:linechargepotential_b}
    \end{subfigure}
        \vspace{-0.6cm}
    \caption{Line segment potential. (a) The variables used in~\cref{eq:linesegmentchargepotential} to compute the potential at point $\boldsymbol{x}$ due to a single layer of charge on a line segment $\Delta S$, parameterized by arclength $s$, from $\boldsymbol{X}_1$ at $s_1$ to $\boldsymbol{X}_2$ at $s_2$. Here $r_1,r_2$ are the distances between $\boldsymbol{x}$ and $\boldsymbol{X}_1,\boldsymbol{X}_2$; $d$ is the perpendicular distance between $\boldsymbol{x}$ and the line segment $\Delta S$; and $\theta_1,\theta_2$ are the angles from the perpendicular line to $\boldsymbol{X}_1,\boldsymbol{X}_2$. (b) The difference between the single layer potential $\mathcal{S}_\Gamma[1]$ for unit charge density on a curve $\Gamma$ and the sum of the charge line potentials $\sum_i \mathcal{S}_{\Delta S_i}[1]$ for 36 segments. The difference is negligible outside the boundary of the ``5$h$ rule''.}
\end{figure}

Following the ``5$h$ rule"~\cite{barnett2014evaluation}, panels within a distance of five times the local grid spacing to an evaluation point are  broken into line segments and evaluated using ~\cref{eq:linesegmentchargepotential}. The contributions to the potential due the other panels are evaluated using ~\cref{eq:uslpsum}. Approximating the curves $\partial \Omega_i$ near the evaluation point as the union of line segments introduces some error especially where the curve is high curvature. This is illustrated in \cref{fig:linechargepotential_b}, which shows the difference between the single layer potential on a curve $\Gamma$ and the approximate potential due to 36 line segments. In practical applications, more line segments are used so that the error introduced by this approximation is much smaller.

\subsection{Exact Solution to Two Non-Concentric Nested Circles Problem}
We will construct a solution to the two non-concentric nested circles problem using a conformal map to a two concentric nested circles problem. The numerical solution to this problem computed using our method and using FEniCS~\cite{alnaes2015FEniCS} will be compared to the exact solution in \cref{sec:testCase1}. The conformal map 
\begin{equation}
\tilde{z} = f(z) = \frac{z - \alpha}{1-\bar{\alpha} z}
\end{equation}
maps the unit disk to the unit disk and the point $\alpha$ to $0$. The inverse of this map is
\begin{equation}
z = f^{-1}(\tilde{z}) = \frac{\tilde{z}+\alpha}{1+\bar{\alpha}\tilde{z}}.
\end{equation}
Let the off-center circle have a diameter on the $x$-axis passing through $(0,0)$ and $(a,0)$. Choose 
\begin{equation}
\alpha = \frac{\sqrt{1+a}-\sqrt{1-a}}{\sqrt{1+a}+\sqrt{1-a}},
\end{equation}
the point equidistant between $(0,0)$ and $(a,0)$ in the unit disc model for hyperbolic geometry.

Let us construct some solutions $\tilde{u}(\tilde{r},\tilde{\theta})$ in the domain with concentric circles. Consider the problem
\begin{align}
    \Delta \tilde{u} &=0~\textrm{on}~\tilde{r}<1, \tilde{r}\neq\alpha\\
    \partial_{\tilde{r}} \tilde{u} &= \sin(m \tilde{\theta})~\textrm{on}~\tilde{r}=1,\\
    \tilde{u}^-&=\tilde{u}^+~\textrm{on}~\tilde{r}=\alpha,\\
    \sigma\partial_{\tilde{r}^-}\tilde{u}&=\partial_{\tilde{r}^+}\tilde{u}~\textrm{on}~\tilde{r}=\alpha,
\end{align}
in which $\sigma=\sigma_2/\sigma_1$ is the ratio of conductivities in the disks of radius $\alpha$ and $1$. $m$ is a natural number. By separating variables on the disk, we find the solution
\begin{equation}
\tilde{u}(\tilde{r},\tilde{\theta}) = \begin{cases} 
      A\tilde{r}^m\sin(m\tilde{\theta}), & \tilde{r} \leq \alpha \\
      (B\tilde{r}^m+C\tilde{r}^{-m})\sin(m\tilde{\theta}), & \alpha \leq \tilde{r} \leq 1
   \end{cases},
   \label{exactSolutionTwoConcentricCircles}
\end{equation}
in which $A = \frac{2}{m} \frac{1}{ \alpha^{2m}(\sigma-1)+\sigma + 1},~~B = A(\sigma+1)/2,~~C = B-1/m$.

We obtain a solution to a non-concentric nested circles problem where $u(r,\theta) = \tilde{u}(\tilde{r},\tilde{\theta})$. The injected current $b_1$ this implies for the non-concentric problem is
\begin{equation}
    b_1(\theta) = \sin\left(m \tan^{-1}\left(\frac{\sin(\theta)(1-\alpha^2)}{\cos(\theta)(1+\alpha^2)-2\alpha}\right)\right) \frac{1-\alpha^2}{1-2\alpha\cos(\theta)+\alpha^2}.
\end{equation}

The charge density on the inner circle is
\begin{equation}
    \gamma_2(\tilde \theta) = \left[\partial_{\tilde{r}^+}\tilde{u} - \partial_{\tilde{r}^-}\tilde{u}\right]_{\tilde{r}=\alpha} = 
\sin(m\tilde\theta)\frac{2\alpha^{m-1}(\sigma-1)}{\alpha^{2m}(\sigma-1)+\sigma + 1}.
\end{equation}

\section{Experimental results}
\label{sec:experiments}

This section illustrates numerically the performance of our method. We compare results obtained using our solver to those obtained using FEniCS~\cite{alnaes2015FEniCS}, a popular open-source (LGPLv3) computing platform for solving partial differential equations. The finite element meshes used in FEniCS were created in Gmsh~\cite{geuzaine2009gmsh}, an open source 3D finite element mesh generator distributed under the terms of the GNU General Public License (GPL). Gmsh allows one to create customizable finite element meshes that can be refined near each interface to increase accuracy at minimal cost to computation time. We consider four separate test cases.

\subsection{Test Case 1: Two Nested Circles}
\label{sec:testCase1}
\begin{figure}[h!]
    \centering
    \begin{subfigure}[b]{0.31\textwidth}
        \centering
        \includegraphics[]{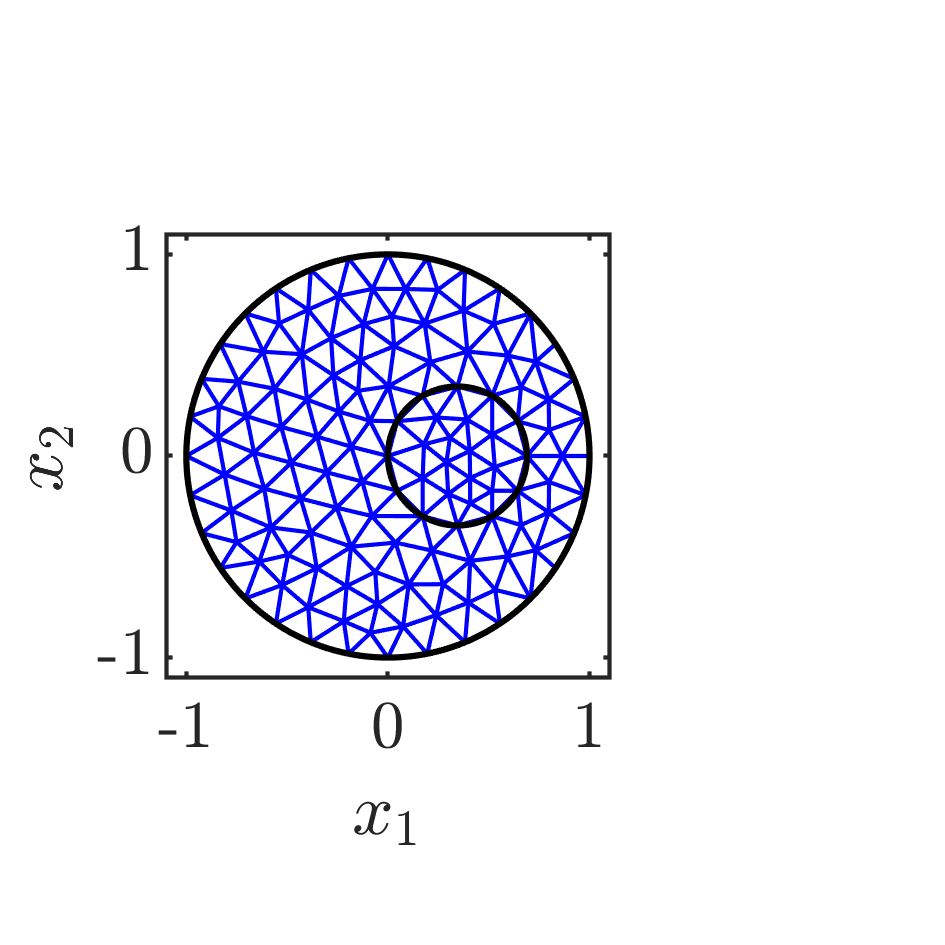}
        \vspace{-1cm}
        \caption{}
    \end{subfigure}
    \hspace{-1.5cm}
    \begin{subfigure}[b]{0.31\textwidth}
        \centering
        \includegraphics[]{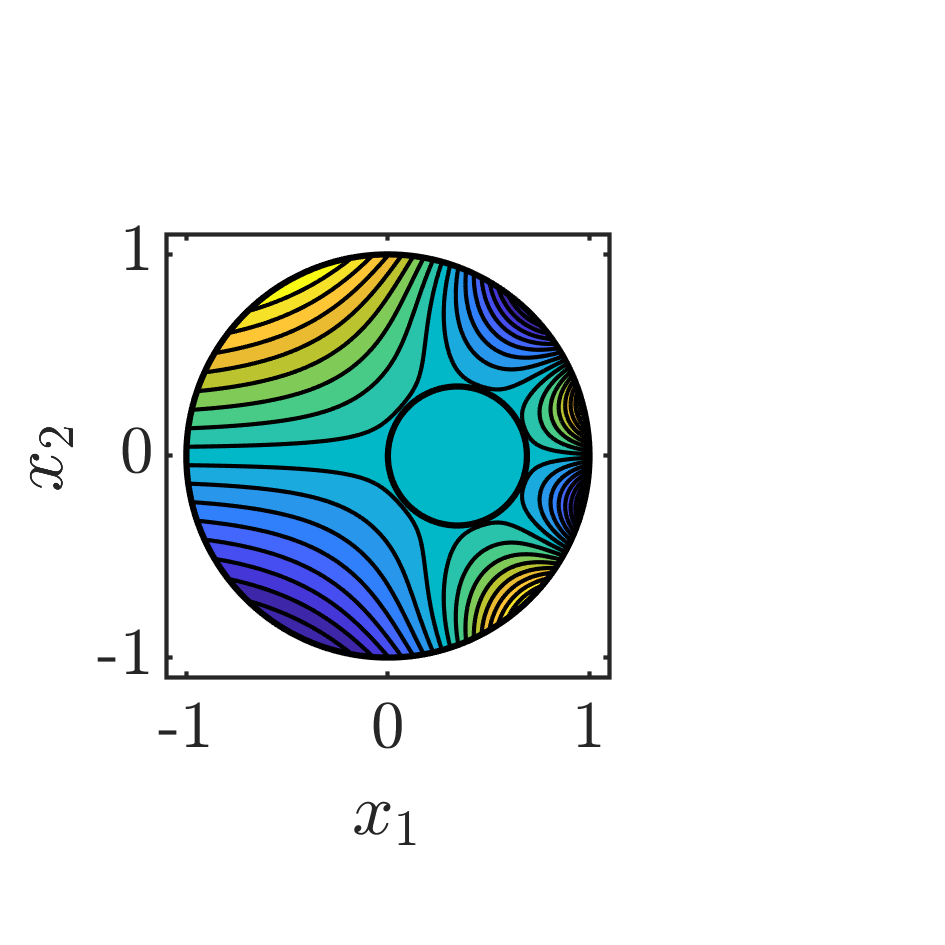}
        \vspace{-1cm}
        \caption{}
    \end{subfigure}
    \hspace{-1.5cm}
    \begin{subfigure}[b]{0.31\textwidth}
        \centering
        \includegraphics[]{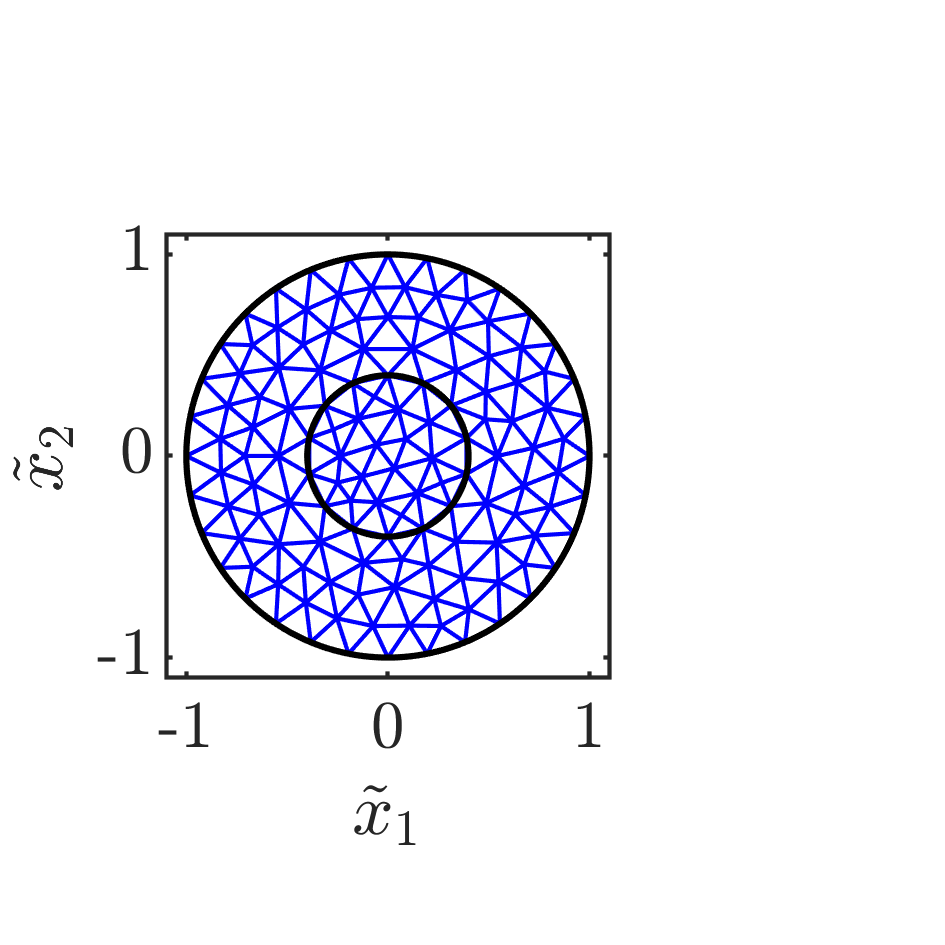}
        \vspace{-1cm}
        \caption{}
    \end{subfigure}
    \hspace{-1.5cm}
    \begin{subfigure}[b]{0.31\textwidth}
        \centering
        \includegraphics[]{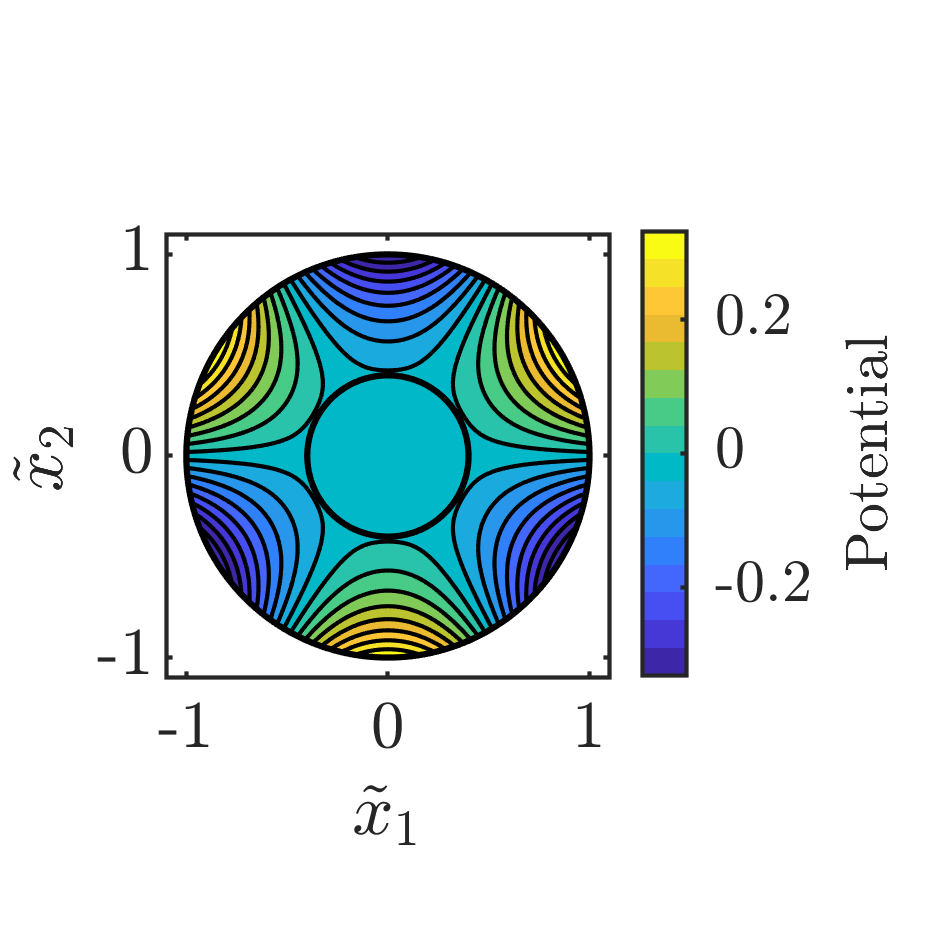}
        \vspace{-1cm}
        \caption{}
    \end{subfigure}\\
    \vspace{-0.4cm}
    \begin{subfigure}[b]{0.31\textwidth}
        \centering
        \includegraphics[]{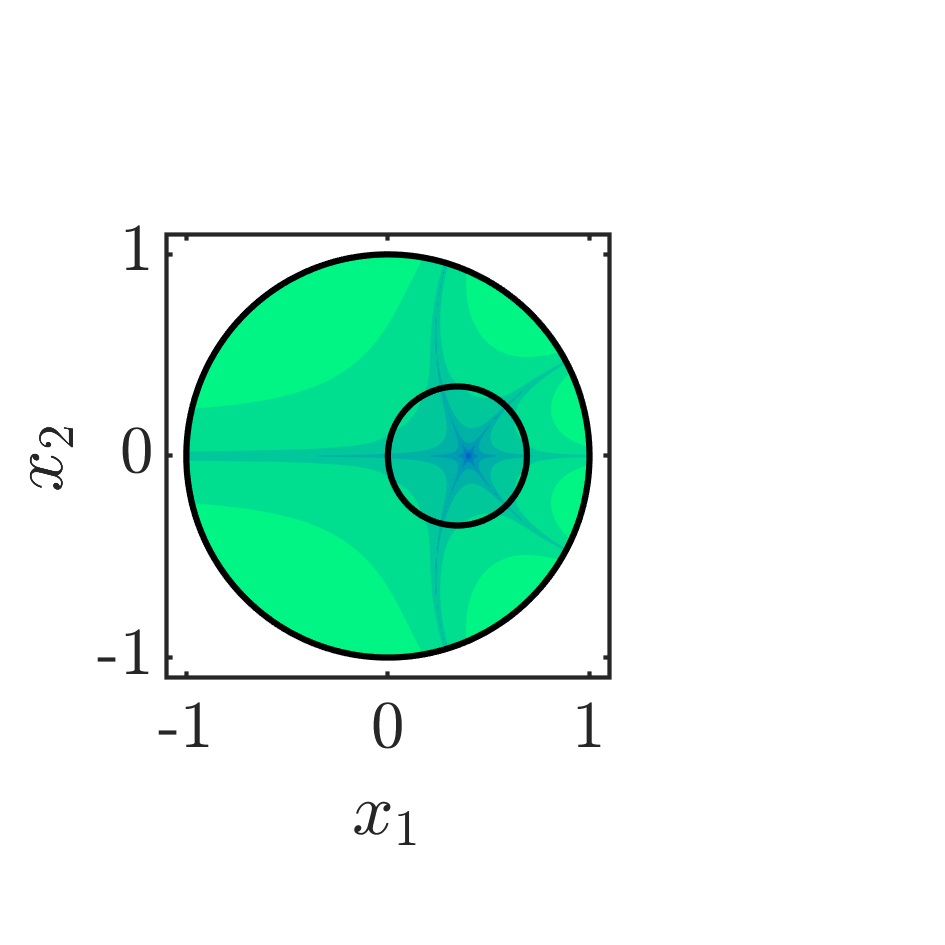}
        \vspace{-1cm}
        \caption{}
    \end{subfigure}
    \hspace{-1.5cm}
    \begin{subfigure}[b]{0.31\textwidth}
        \centering
        \includegraphics[]{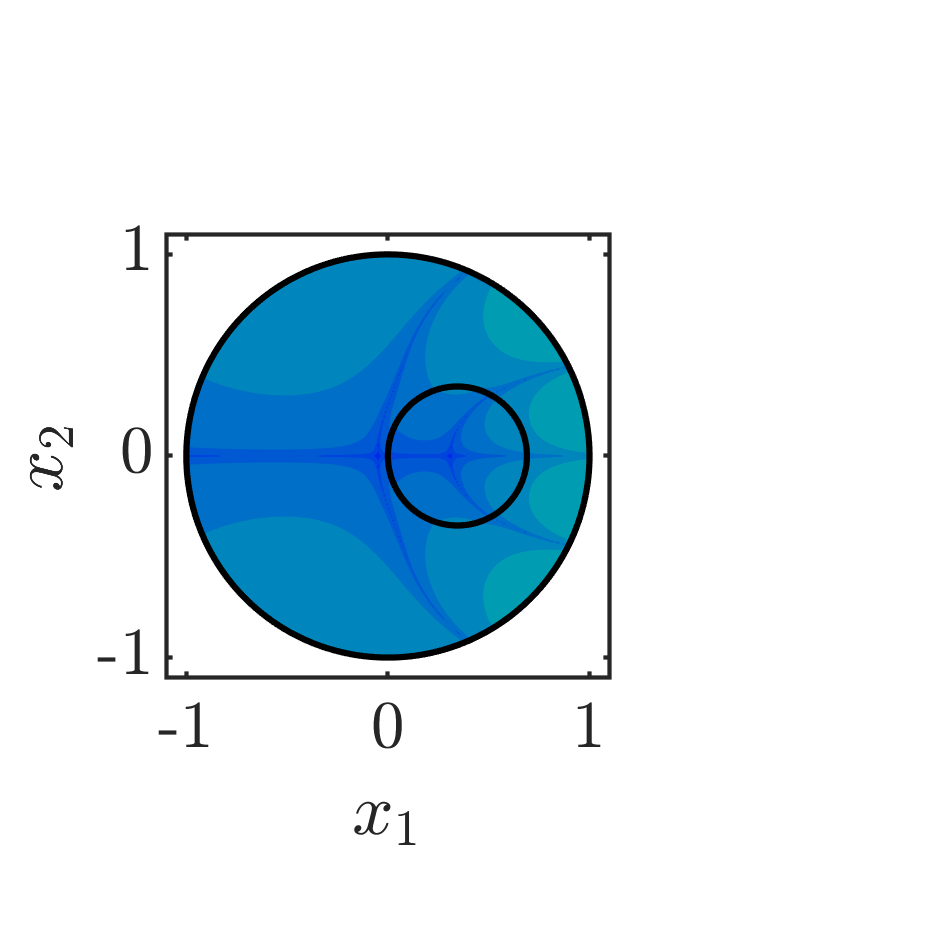}
        \vspace{-1cm}
        \caption{}
    \end{subfigure}
    \hspace{-1.5cm}
    \begin{subfigure}[b]{0.31\textwidth}
        \centering
        \includegraphics[]{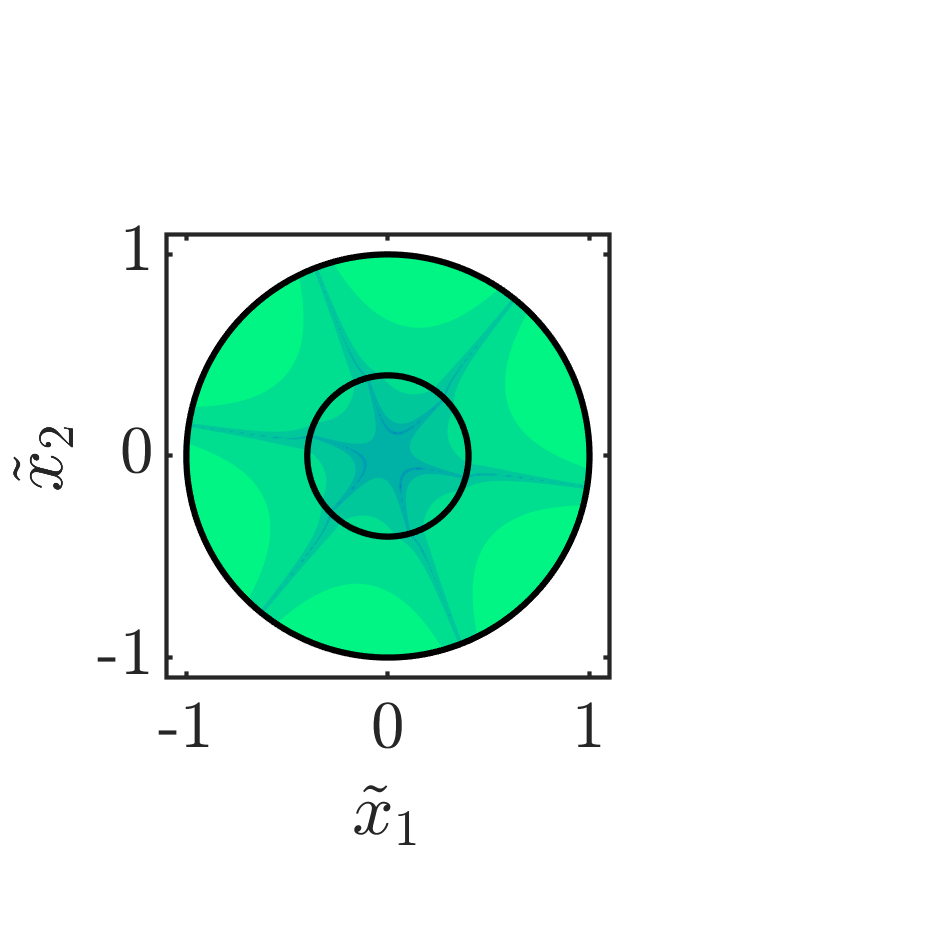}
        \vspace{-1cm}
        \caption{}
    \end{subfigure}
    \hspace{-1.5cm}
    \begin{subfigure}[b]{0.31\textwidth}
        \centering
        \includegraphics[]{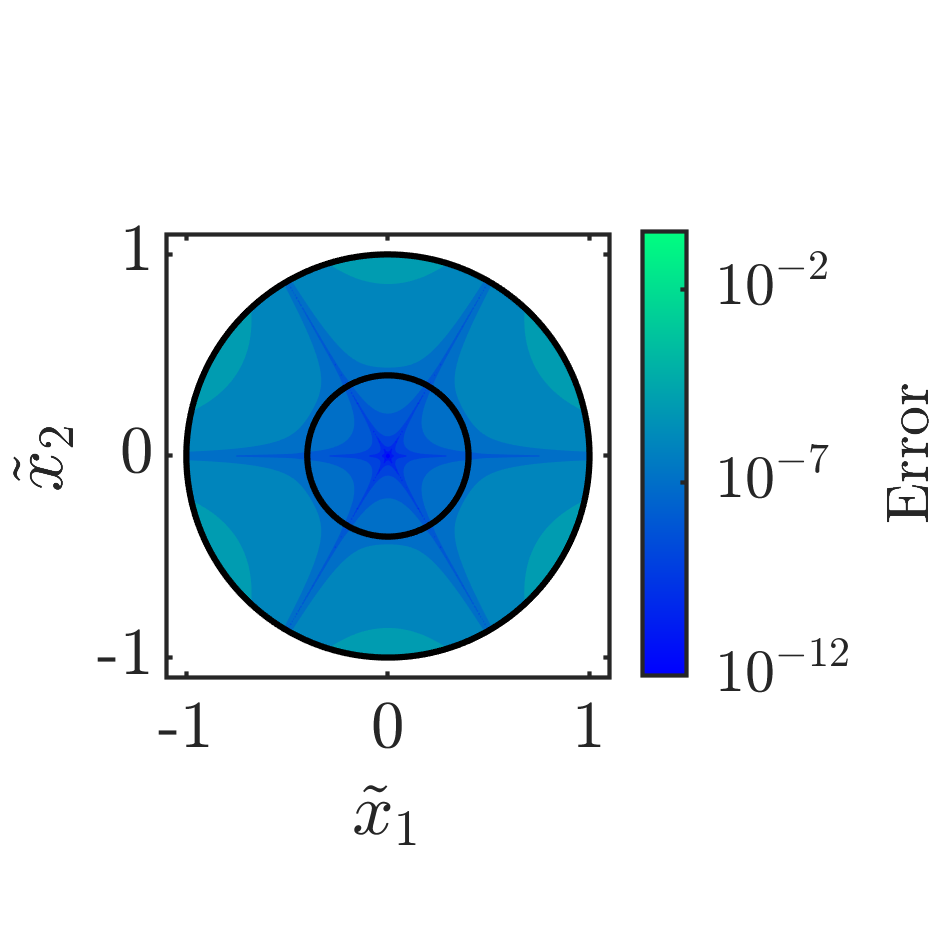}
        \vspace{-1cm}
        \caption{}
    \end{subfigure}
    \caption{Pointwise errors for the nested circles test case. In the non-concentric case, a coarse version of the finite element mesh created in Gmsh and used by FEniCS (a), the contours of the solutions (b), the pointwise error of the numerically computed solutions using FEniCS (e) and using our method (f). The corresponding plots for the concentric case are shown in (c,d,g, and h) respectively. Parameter values are $m=3,\sigma=2,\alpha=0.40$ ($a\approx0.69$).
    }
    \label{fig:twocirclestest}
\end{figure}

\begin{figure}[h!]
     \centering
     \begin{subfigure}[b]{0.49\textwidth}
         \centering
         \includegraphics[]{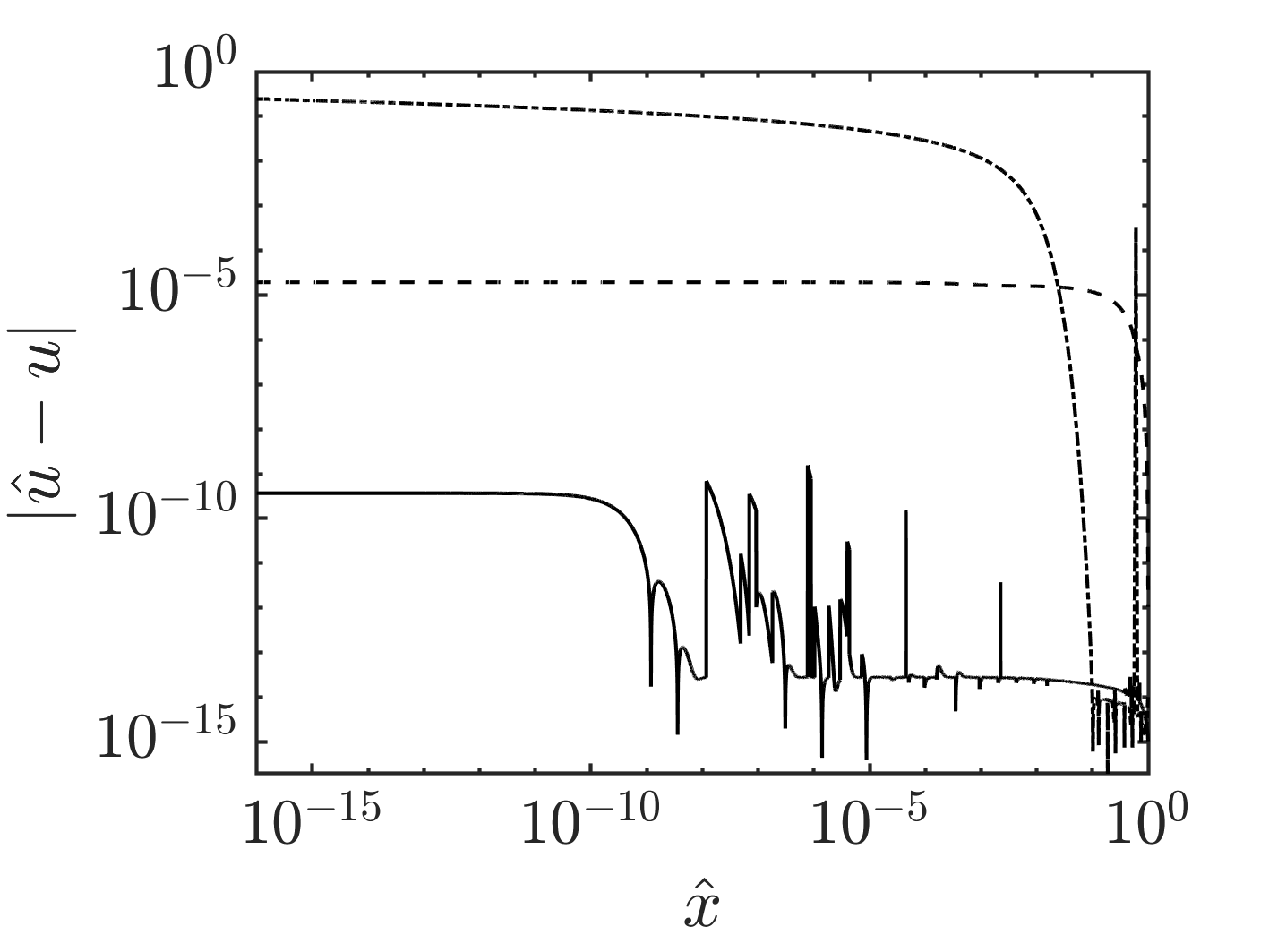}
        \vspace{-0.2cm}
         \caption{}
         \label{fig:case1evaluation_atnode}
     \end{subfigure}
     \hfill
     \begin{subfigure}[b]{0.49\textwidth}
         \centering
         \includegraphics[]{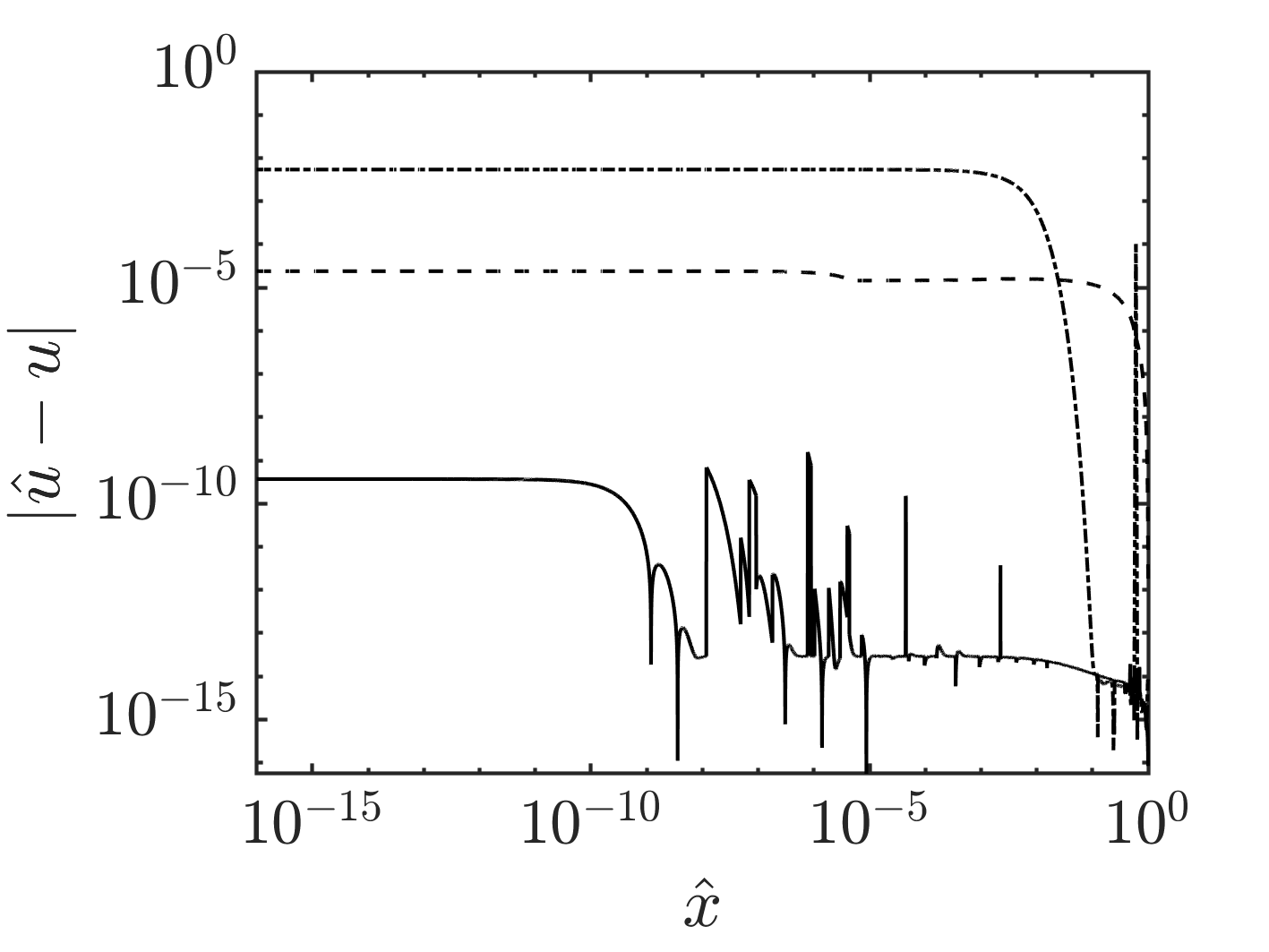}
        \vspace{-0.2cm}
         \caption{}
         \label{fig:case1evaluation_betweennode}
     \end{subfigure}
     \\
        \vspace{-0.6cm}
        \caption{Absolute error in the computed solution using the na\"{i}ve quadrature (dash-dotted), using adaptive quadrature using the interpolant of the charge densities (solid), and using the charge line potential between each node (dashed) for the case of two concentric circles as a function of distance from the outer boundary. The solution is evaluated at the point $(r,\theta)=(1-\hat{x},\frac{\pi}{2})$ in polar coordinates. In (a) the point $(r,\theta)=(1,\frac{\pi}{2})$ is at a quadrature node and in (b) we rotate the quadrature nodes so that the same point is half way between two quadrature nodes. On each interface/boundary $2^8$ uniformly spaced nodes were used. For the adaptive quadrature, an adaptive tolerance of $10^{-6}$ was used.}
        \label{fig:case1evaluation}
\end{figure}

We consider the case of two nested circles. \Cref{fig:twocirclestest} shows the analytic solution to a two circle non-concentric and the corresponding concentric problems, the pointwise difference between the analytic solutions and the boundary integral method solutions, and the pointwise difference between the analytic solutions and the solutions computed using FEniCS. Additionally, the  finite element mesh created using the frontal Delaunay triangulation~\cite{rebay1993efficient} in Gmsh and used by FEniCS is shown. The wall time to build and solve the system of equations for the two non-concentric circle  problem was 0.23 seconds for our solver and 1.37 seconds for FEniCS and the two concentric circle problem was 0.07 seconds for our solver and 0.83 seconds for FEniCS. \Cref{fig:case1evaluation} shows the absolute error in the computed solution using the na\"{i}ve quadrature (dash-dotted), using adaptive quadrature using the interpolant of the charge densities (solid), and using the charge line potential between each node (dashed) for the case of two concentric circles as a function of distance from the outer boundary. The solution is evaluated at the point $(r,\theta)=(1-\hat{x},\frac{\pi}{2})$ in polar coordinates. In \cref{fig:case1evaluation_atnode} the point $(r,\theta)=(1,\frac{\pi}{2})$ is at a quadrature node and in \cref{fig:case1evaluation_betweennode} we use the same number of quadrature nodes but rotate them so that the same point $(r,\theta)=(1,\frac{\pi}{2})$ is half way between two quadrature nodes. The wall time to evaluate the solution at all 8,000 points used in \Cref{fig:case1evaluation} was 0.02 seconds when using the na\"{i}ve quadrature, 36.16 seconds when using adaptive quadrature, and 0.90 seconds when using the charge line potential. The charge line potential leads to more accurate results near the boundary than the na\"{i}ve quadrature and is faster than adaptive quadrature. The case of two nested circles is a standard test problem that has been solved previously via finite element methods~\cite{gehre2014analysis}.

\subsection{Test Cases 2, 3, and 4: Examples of General Curves}
\begin{figure}
    \centering
        \hspace{-0.75cm}
    \begin{subfigure}[b]{0.43\textwidth}
        \centering
        \includegraphics[]{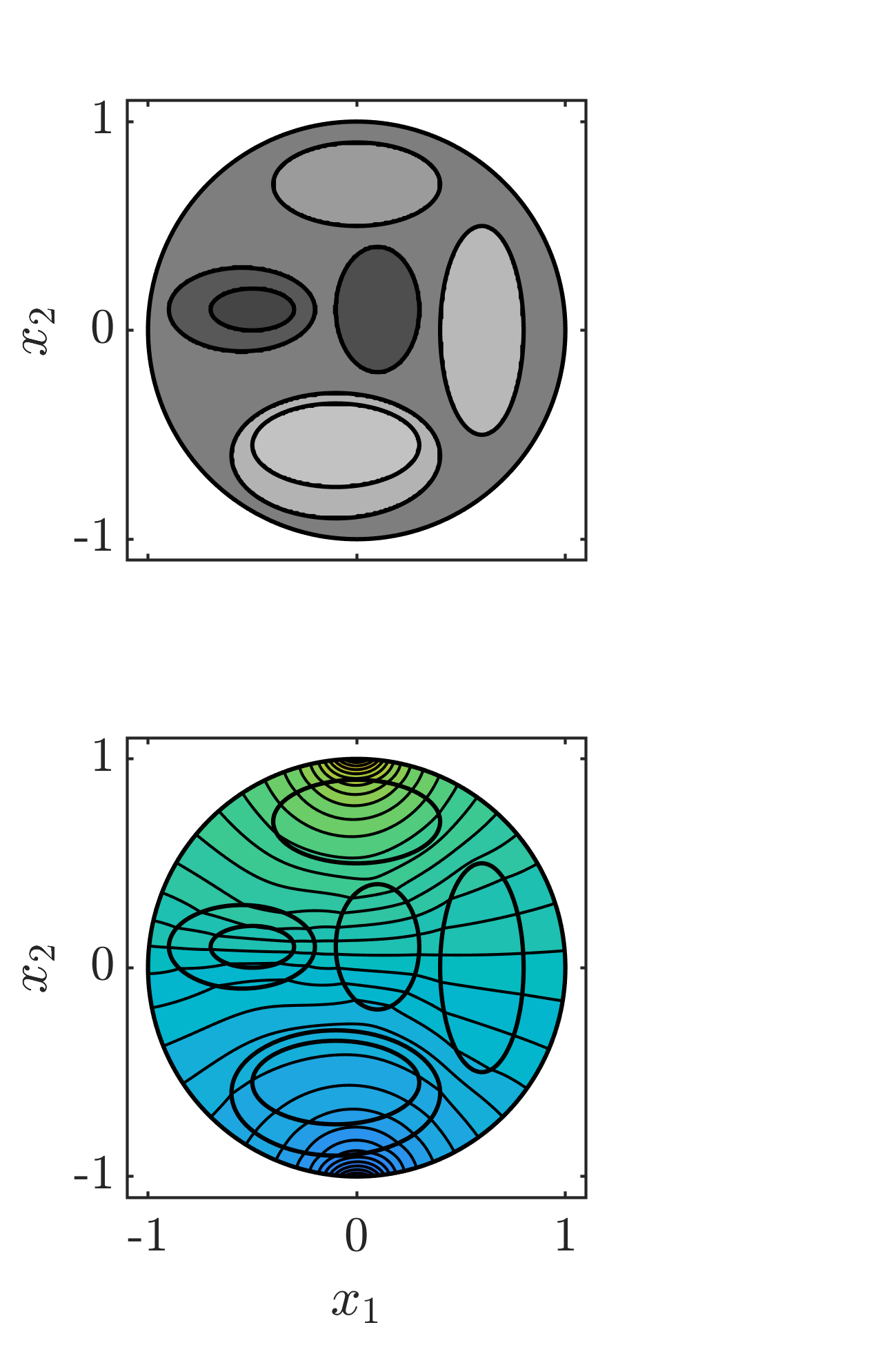}
        \vspace{-0.7cm}
        \caption{}
        \label{fig:sigma_potential_2}
    \end{subfigure}
        \hspace{-2.2cm}
    \begin{subfigure}[b]{0.43\textwidth}
        \centering
        \includegraphics[]{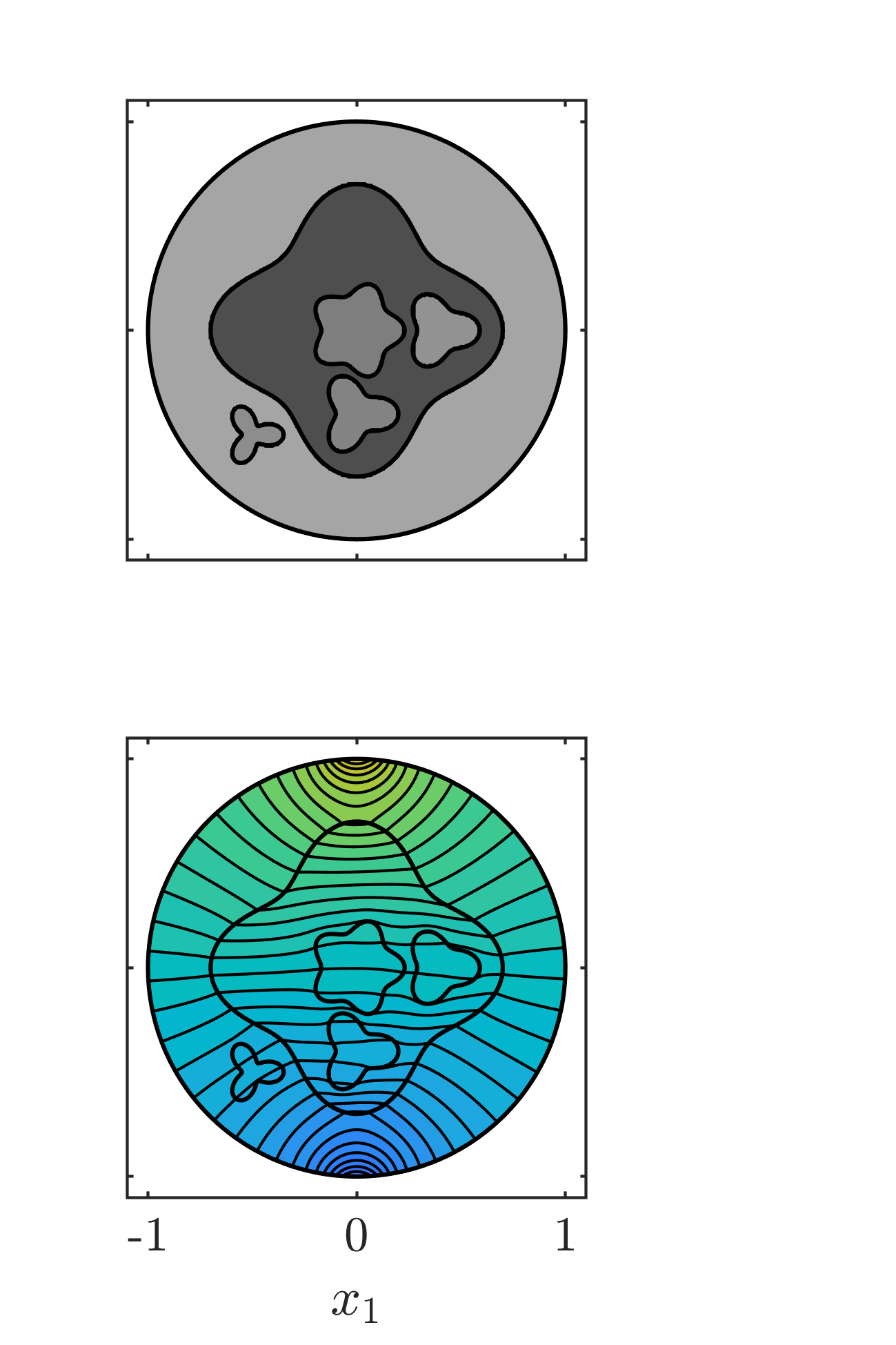}
        \vspace{-0.7cm}
        \caption{}
        \label{fig:sigma_potential_3}
    \end{subfigure}
        \hspace{-2.2cm}
    \begin{subfigure}[b]{0.43\textwidth}
        \centering
        \includegraphics[]{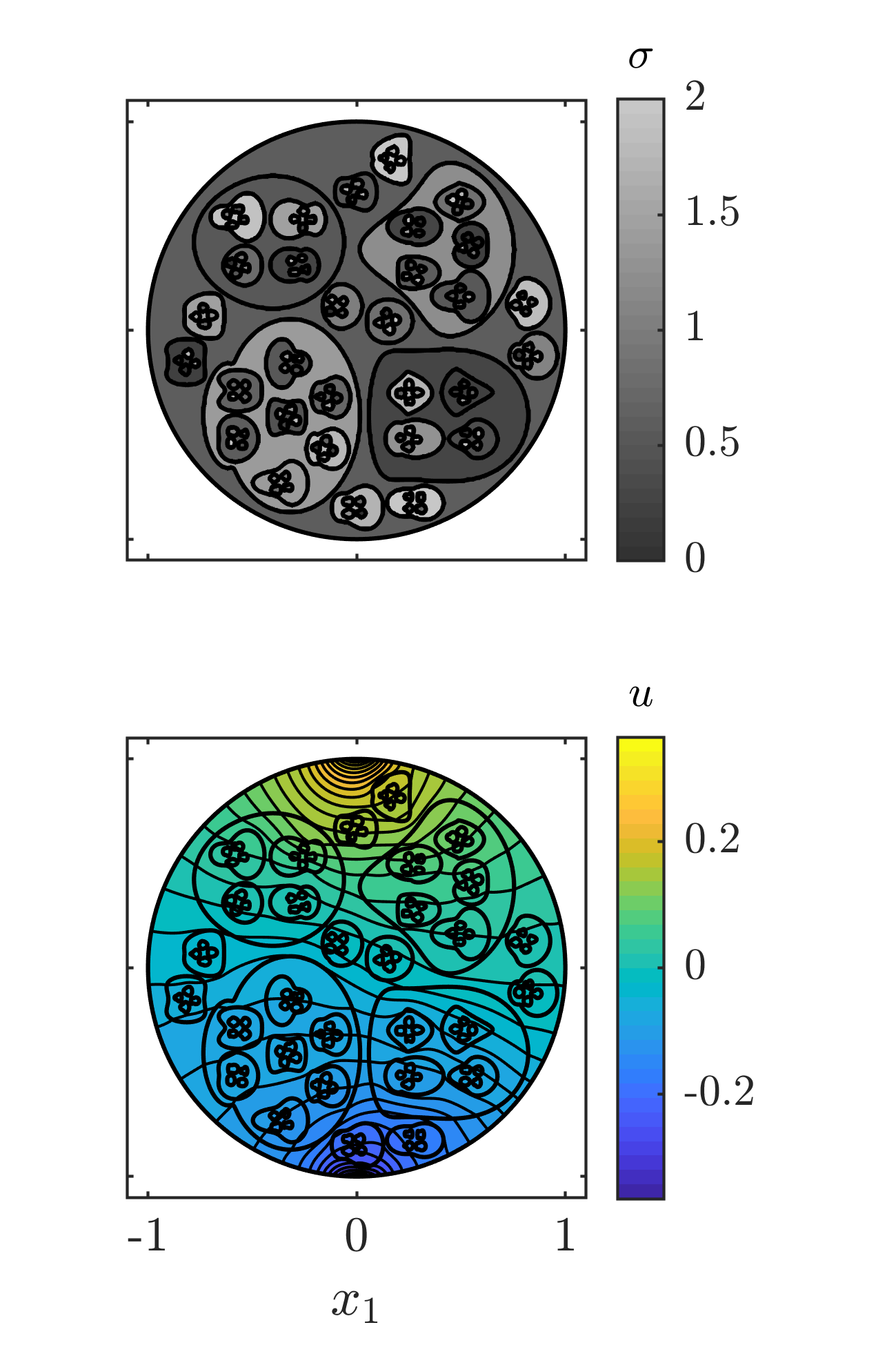}
        \vspace{-0.7cm}
        \caption{}
        \label{fig:sigma_potential_4}
    \end{subfigure}
    \\
    \vspace{-0.6cm}
    \caption{Conductivity $\sigma$ for the problem (top) and potential $u$ computed using our solver (bottom) for test case 2 (a), test case 3 (b), and test case 4 (c). In all three test cases the injected current on the outer boundary is $b_1(\theta)=\cos(6\theta-\pi)\1_{\{|\theta-\pi/2|<\pi/12\}}-\cos(6\theta-\pi)\1_{\{|\theta-3\pi/2|<\pi/12\}}$.}
    \label{fig:test_cases_sigma_potential}
\end{figure}

We consider cases where the conductivity is piecewise constant on regions defined by more general curves. Test case 2 has elliptical regions, test case 3 has clover leaf/starfish curves where the radius of each curve is given in polar coordinates by $r(\theta) = A + B \cos(C\theta)$ for various parameters $A,B,C$, and test case 4 has 155 different regions of constant conductivity. In each case, the domain is the unit disk and the injected current on the outer boundary is $b_1(\theta)=\cos(6\theta-\pi)\1_{\{|\theta-\pi/2|<\pi/12\}}-\cos(6\theta-\pi)\1_{\{|\theta-3\pi/2|<\pi/12\}}$. \Cref{fig:test_cases_sigma_potential} displays the regions of constant conductivity and the computed potential using our solver for test cases 2, 3, and 4. A similar problem to test case 2 involving elliptical regions of constant conductivity has been solved using an integral equation method~\cite{nasser2011boundary}.

\subsection{Grid Refinement and Rescaling Studies}
\begin{figure}
     \centering
     \begin{subfigure}[b]{0.33\textwidth}
         \centering
         \includegraphics[]{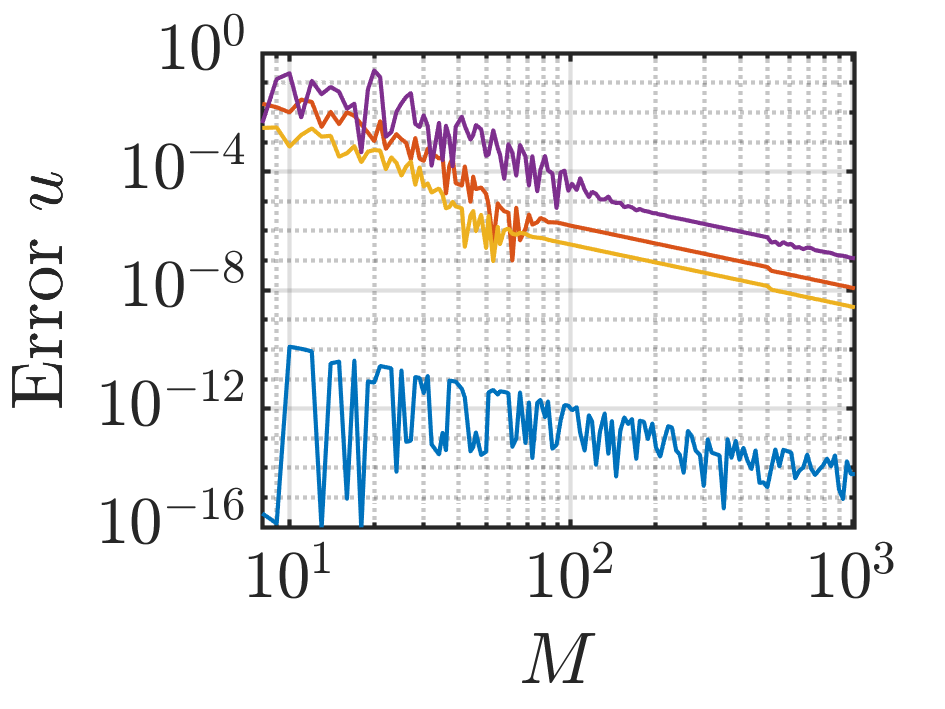}
        \vspace{-0.2cm}
         \caption{}
         \label{fig:testCase1234_refinement_u}
     \end{subfigure}
        \hspace{-0.2cm}
     \begin{subfigure}[b]{0.33\textwidth}
         \centering
         \includegraphics[]{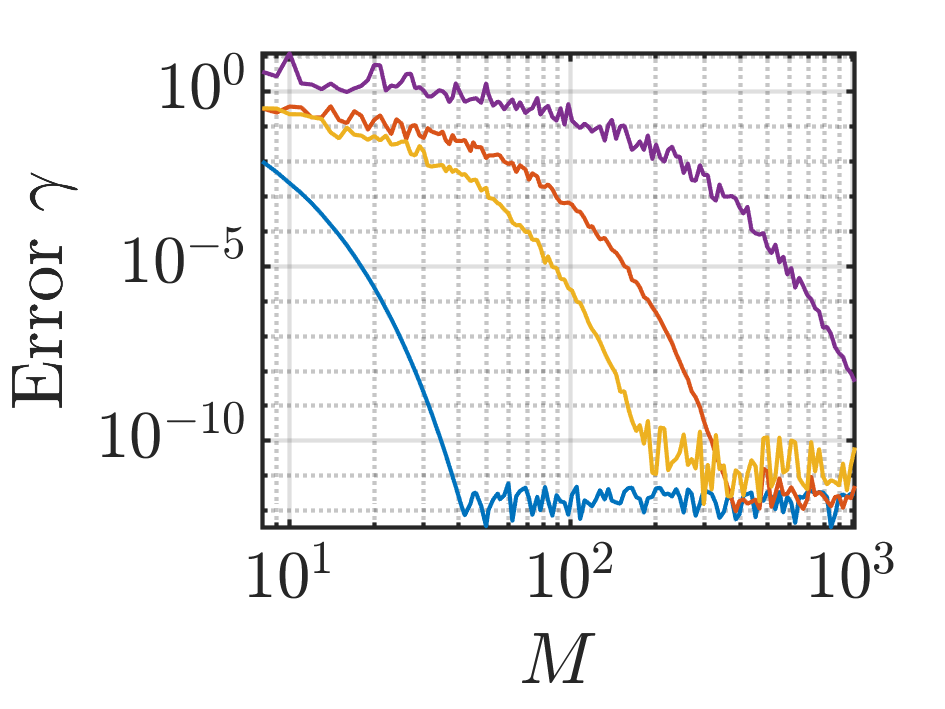}
        \vspace{-0.2cm}
         \caption{}
         \label{fig:testCase1234_refinement_gamma}
     \end{subfigure}
        \hspace{-0.2cm}
     \begin{subfigure}[b]{0.33\textwidth}
         \centering
         \includegraphics[]{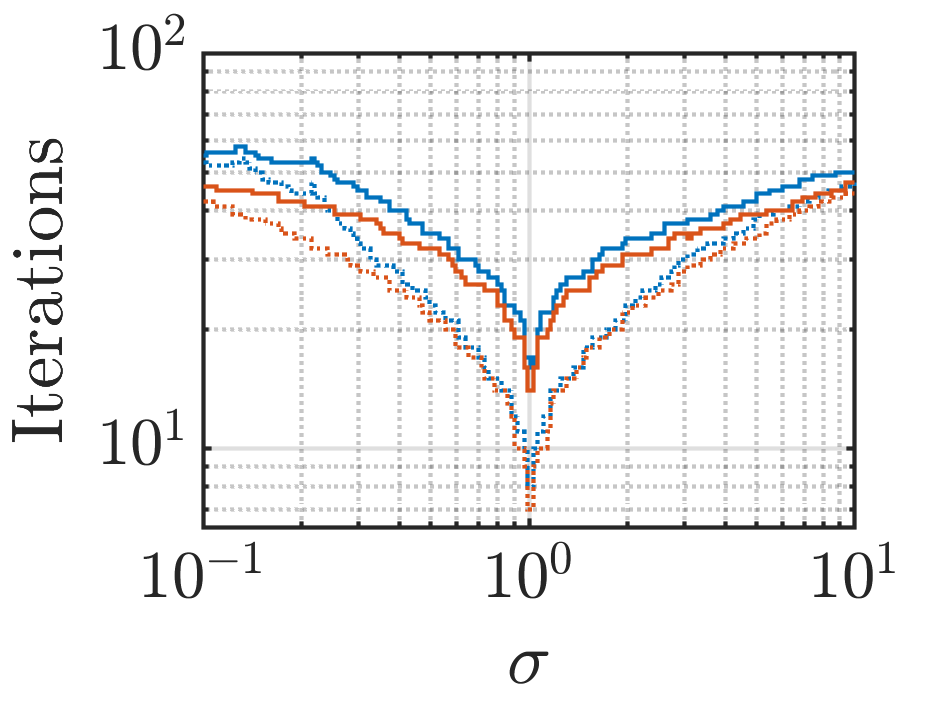}
        \vspace{-0.2cm}
         \caption{}
         \label{fig:testCase2_refinement_scaling}
     \end{subfigure}
     \\
        \vspace{-0.6cm}
        \caption{Refinement studies showing (a) the pointwise error estimate of the solution $u$ at the boundary point $(1,0)$ and (b) the error estimate of the charge density $\gamma$ as a function of $M$, the number of uniform grid points on each interface, for all four test cases -- 1 (blue), 2 (red), 3 (yellow), 4 (purple). (c) Number of GMRES iterations (with GMRES tolerance, \textit{i}.\textit{e}., relative residual error, of $10^{-8}$ and initial guess of the zero vector) performed with (dotted lines) and without (solid lines) rescaling the charge densities for seven elliptical regions of constant conductivity alternating between 1 and $\sigma$. On each of the seven ellipses, $M=2^5$ (blue lines) or $2^8$ (red lines) uniformly spaced nodes were used.}
        \label{fig:refinement_and_scaling}
\end{figure}

\Cref{fig:testCase1234_refinement_u,fig:testCase1234_refinement_gamma} display two refinement studies for our solver. In the first refinement study, the error estimate of the potential $u$ is plotted as a function of $M$, the number of uniform grid points on each interface, for each of the four test cases. The error estimate is taken to be the absolute difference between the potential difference between the solutions computed at the boundary point $(1,0)$ with $M$ and $2M$ uniform grid points on each interface. The error is primarily due to solution evaluation error. If higher accuracy in the evaluation is required then using adaptive quadrature will improve the accuracy of the solver albeit at an increased cost to computation time. In the second refinement study, the charge density at each node on each interface/boundary is computed for each $M$ and the maximal difference between these charge densities and the charge densities computed with $2M$ grid points is plotted. The refinement studies show that the charge densities converge spectrally and the potential converges linearly for our solver for all four test cases.

As discussed in \cref{sec:Conditioning of Linear System}, especially in the case of close conductivities, scaling the charge densities improves the conditioning of the resulting linear system. A study on the effect of rescaling the charge densities for our solver was performed for the case of the seven elliptical regions of constant conductivity described in test case 2 except with the conductivities alternating between 1 and a variable value $\sigma$ for nested regions. The number of GMRES iterations performed with and without rescaling the charge densities as a function of the value $\sigma$ are shown in \cref{fig:testCase2_refinement_scaling}. Rescaling decreased the number of GMRES iterations required. Note that iterations with a larger number of nodes per interface take longer to complete.

\begin{figure}
     \centering
     \begin{subfigure}[b]{0.49\textwidth}
         \centering
         \includegraphics[]{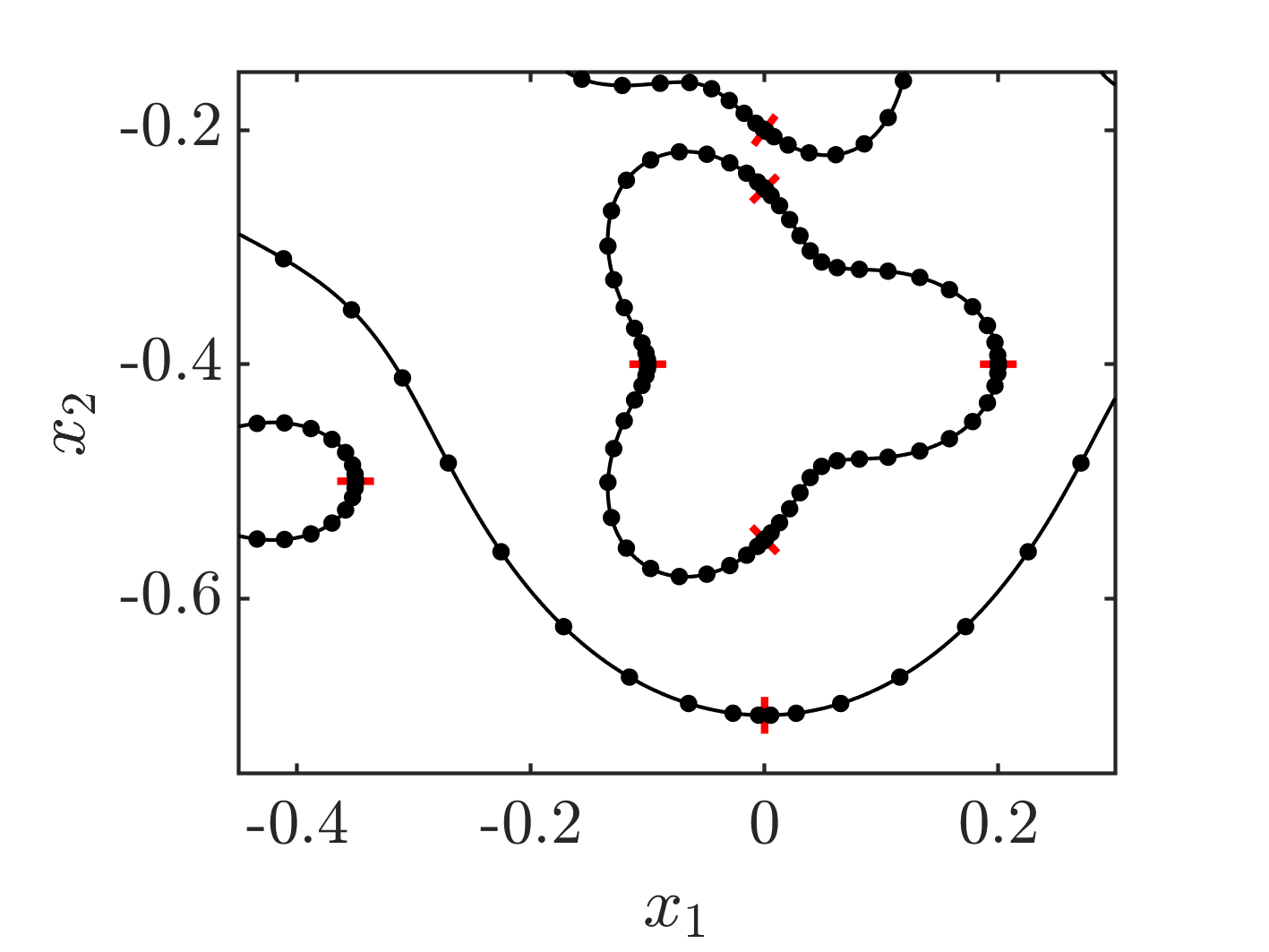}
        \vspace{-0.2cm}
         \caption{}
         \label{fig:adaptive_grid}
     \end{subfigure}
     \hfill
     \begin{subfigure}[b]{0.49\textwidth}
         \centering
         \includegraphics[]{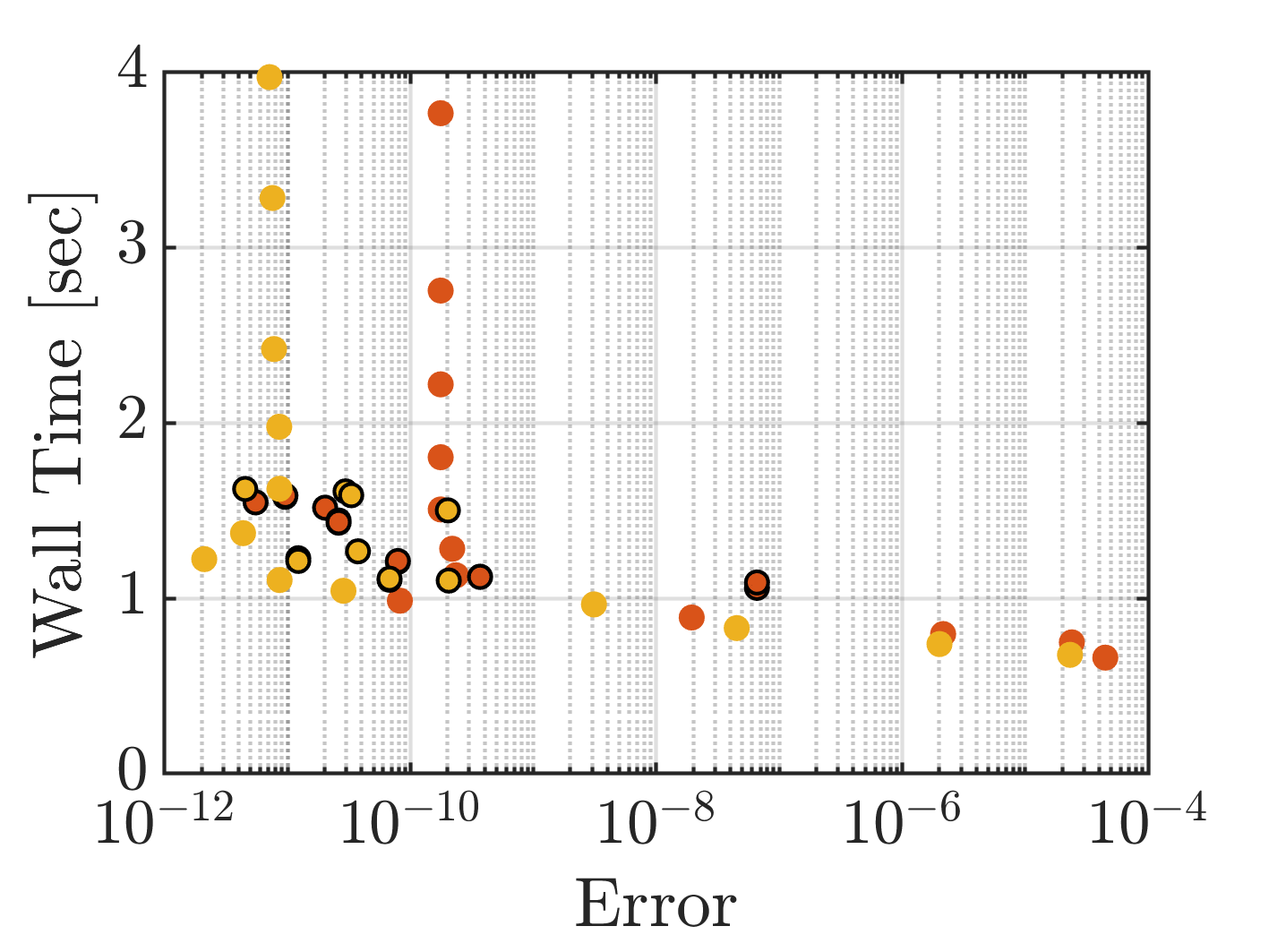}
        \vspace{-0.2cm}
         \caption{}
         \label{fig:refinement_performance_plot}
     \end{subfigure}
     \\
        \vspace{-0.6cm}
        \caption{(a) Part of the automatically selected adaptive grid (black dots are the quadrature nodes, red lines are the panel boundaries, black curves are interfaces) for test case 3. (b) Refinement performance plot showing wall time averaged over 10 samples vs error estimate for our solver with (black outline) and without (no outline) adaptively selecting the grid points for test case 2 (red) and 3 (yellow).}
        \label{fig:adaptive_grid_refinement_performance_plot}
\end{figure}

In order to benefit from the improved accuracy that using a larger number of nodes per interface brings without suffering too much increased computation time, we adaptively refine our grid to place more points where the interfaces are close together. \Cref{fig:adaptive_grid} shows an example of the automatically selected adaptive grid for test case 3. \Cref{fig:refinement_performance_plot} demonstrates how the adaptive grid refinement improves accuracy at minimal cost to computation time for test cases 2 and 3. The total wall time for the build and solve averaged over 10 samples is plotted against an error estimate -- the absolute difference between the potential computed at the origin for the given run and a run computed with higher accuracy. For both test cases, the adaptive method was more efficient ultimately achieving more accurate results with less computation time. All simulations were performed on a desktop PC with a 6-Core 3.59 GHz CPU and 16 GB RAM.

\section{Existence and Regularity}
\label{sec:Existence and Regularity}

In this section, we will provide some theoretical study of the system of integral equations \crefrange{eq:inteqn1}{eq:inteqn2}. The regularity results we obtain on the charge densities guide our choice of quadrature for discretizing the integral equations. Since our adaptive method is based on approximating the charge densities by their truncated Fourier series or truncated Legendre series, the guarantees on the smoothness of the charge densities that we derive ensure the series approximations converge quickly. The discretization scheme and details of the solution evaluation were discussed in depth in \cref{sec:methods}. Some conclusions follow in \cref{sec:conclusions}. For now, we aim to show two separate results.

First, we will establish existence and uniqueness to this system of equations. This will be obtained from general PDE theory, and the analogous statement for the system \cref{eq:inteqn1,eq:inteqn2} will follow as a consequence. 

 Second, we will derive regularity for the charge densities along each interface. The latter result will have two formulations. In both, we will obtain regularity of the charge density functions along an interface, assuming the interface itself enjoys the same regularity. The most general result will be \cref{thm:cgeneral} below, which will show that assuming that the interface $\partial\Omega_i$ enjoys $\mathcal{C}^k$ regularity, then the charge density $\gamma_i(\boldsymbol{x})$ is of regularity $H^{k}$ (i.e. in the Hilbert space of order $k$). In particular, it will also be in the H\"{o}lder space $\mathcal{C}^{k-2,\alpha}$, for $0\le \alpha<1$ by Sobolev embedding. In this theorem, the control we obtain on these norms will be in terms of the $L^2$-norm of the injected current. In the supplementary materials, we present \cref{sm:thm:concentric} where we assume the regions of constant conductivity are concentric circles and obtain stronger control on the charge densities.

\begin{theorem}
\label{thm:cgeneral}
Consider the solution of problem \cref{eq:inteqn1,eq:inteqn2}, in which $\partial\Omega_i$ are general ${\cal C}^k$ curves and $b_1\in L^2(\partial\Omega)$. Consider the charge densities $\gamma_i(\cdot )$ at each of the inner interfaces. Then these functions are ${\cal C}^{k-2,\alpha}$ regular for $0\le\alpha<1$, and their norms are bounded as follows: ($s$ is the arc-length parameter on each curve)
\begin{equation}
\|\gamma_i(s)\|_{{\cal C}^{k-2,\alpha}(\partial\Omega_i)}\le \left|1-\frac{\sigma_{p_i}}{\sigma_i}\right| \frac{M_i}{\min\{\sigma_1,\dots,\sigma_N\}} \rho_{\sigma_i}\|b_1\|_{ L^2(\partial\Omega)}
\end{equation}
Here the constant $M_i$ depends on: The shape of the interface $\partial\Omega_i$, the distance $\delta_i$ between the curve $\partial\Omega_i$ and its nearest neighboring curves, and the order $k\in\mathbbm{N}$. $\rho_{\sigma_i}$ is the maximum ratio of neighboring conductivities \textit{i}.\textit{e}., $\rho_{\sigma_i}:=\textrm{max}\{\frac{\sigma_i}{\sigma_{p_i}},\frac{\sigma_{p_i}}{\sigma_{i}}\}$.
\end{theorem}
In particular, when an interface $\partial\Omega_i$  is very close to another interface $\partial\Omega_j$, ($j\ne i$) the constants $M_i$ will increase. The precise constant $M_i$ can be obtained by following the proof below.

\subsection[Regularity of the solution at the interfaces, in the context of general smooth interfaces]{Regularity of the solution at the interfaces, in the context of general ${\cal C}^k$ interfaces}

Now we examine regularity of the solution for general ${\cal C}^k$ interfaces. Let $b_1\in L^2(\partial\Omega)$. The weak formulation of \cref{eq: PDE,eq: BC} is
\begin{equation}
\label{eq: weak formulation general curves}
    \int_\Omega \sigma \langle \nabla u, \nabla \phi \rangle~dA = \int_{\partial \Omega} b_1 \phi~dl
\end{equation}
for all $\phi\in H^{1}(\Omega)$. Since $\sigma\in L^\infty$, from general PDE theory we obtain that there exists a unique weak solution $u\in H^{1}_0(\Omega)$~\cite{evans10}. Setting $\phi=u$ and using the Cauchy-Schwarz inequality, the trace theorem, and the Poincar\'{e} inequality as in the proof of \cref{sm:eq:bound_on_grad_u} in the supplementary material, we obtain
\begin{align}\label{eq:bound_on_grad_u_general}
    \|\nabla u\|_{L^2(\Omega)} &\le \frac{c}{\min\{\sigma_1,\dots,\sigma_N\}}\|b_1\|_{L^2(\partial\Omega)}
\end{align}
in which $c$ is a constant depending on the trace constant and Poincar\'{e} constant for the domain $\Omega$.

First, we will examine tangential regularity. 
Let $w\in C^\infty(\Omega)$ and let $T$ be a vector field that on $\partial\Omega_i$ is of unit length and tangent to $\partial\Omega$ and is zero before the next curve. (Any such vector field $T$ will do for now --- we further fix $T$ momentarily). 
Then, we take $\phi=T(w)$ in \cref{eq: weak formulation general curves} \textit{i}.\textit{e}.,
\begin{equation}
        \int_\Omega \sigma \langle \nabla u, \nabla T (w) \rangle~dA = \int_{\partial \Omega} b_1 T (w)~dl=0.
\end{equation}
We integrate by parts with respect to $T$ to obtain: 
\begin{align}
        \int_\Omega \sigma \langle \nabla T(u), \nabla w \rangle~dA &= - \int_{\Omega}\sigma
    \cdot\textrm{div}(T)\cdot\langle \nabla u,\nabla w\rangle ~dA
    +2\int_{\Omega} \sigma \nabla_aT_b \nabla^au \nabla^bw ~dA.
\end{align}
(We use the Einstein summation convention so $a$ and $b$ are summer from 1 to 2). 

By formally setting $w=T(u)$ and using the Cauchy–Schwarz inequality, we obtain the inequality: 
\begin{align}
  {\rm min}|\sigma|\cdot \sqrt{  \int_{\Omega}|\nabla T(u)|^2~dA}&\le {\rm max}|\sigma|\cdot\left(\sup|{\textrm{div}(T)}|+2\sup|\nabla T|\right)\sqrt{\int_{\Omega}|\nabla u|^2~dA}.
\end{align}

Thus $\nabla T (u)\in L^2(\Omega)$, with the bound: 
\begin{equation}
    \|\nabla T(u)\|_{L^2(\Omega)}\le
    4\sup|\nabla T|\cdot\frac{{\rm max}|\sigma|}{{\rm min}|\sigma|}\|\nabla u\|_{L^2(\Omega)}
\end{equation}
Iterating, we get that the composition
\begin{equation}
    \label{eq:higher norm}
T^j(u)\in H^1(\Omega)
\end{equation}
for all $j\le k$. The constant for each such iteration will increase by a factor 
\begin{equation}
    4\sup|\nabla T|\cdot\frac{{\rm max}|\sigma|}{{\rm min}|\sigma|}
\end{equation}
\textit{i}.\textit{e}.,
\begin{equation}\label{eq:bound_nabla_T_j_u}
    \|\nabla T^j(u)\|_{L^2(\Omega)}\le
    \left(4\sup|\nabla T|\cdot\frac{{\rm max}|\sigma|}{{\rm min}|\sigma|}\right)^j\|\nabla u\|_{L^2(\Omega)}
\end{equation}
for $j\le k$. In particular, we have obtained extra regularity in the direction of the vector field $T$. 

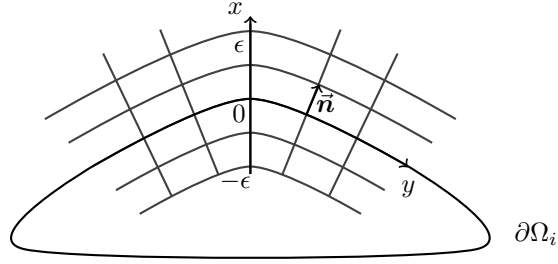
\begin{figure}[h!]
    \centering
    \begin{tikzpicture}
     \foreach \x in {0.7,0.85,1.0,1.15,1.3}
        \draw [darkgray, thick] plot [smooth] coordinates {(-2.1*\x,2.1*\x) (0,3*\x) (2.1*\x,2.1*\x)};
    \tikzmath{\r1 = 2; \r2=4.1; \r3=2.955;} 
    \draw [darkgray, thick] plot [smooth] coordinates {({2*cos(58.5)},{2*sin(58.5)}) ({4.1*cos(61.5)},{4.1*sin(61.5)})};
    \draw [darkgray, thick] plot [smooth] coordinates {({2.03*cos(78)},{2.03*sin(78)}) ({4.1*cos(73)},{4.1*sin(73)})};
    \draw [darkgray, thick] plot [smooth] coordinates {({2*cos(90)},{2*sin(90)}) ({4*cos(90)},{4*sin(90)})};
    \draw [darkgray, thick] plot [smooth] coordinates {({2.03*cos(102)},{2.03*sin(102)}) ({4.1*cos(107)},{4.1*sin(107)})};
    \draw [darkgray, thick] plot [smooth] coordinates {({2*cos(121.5)},{2*sin(121.5)}) ({4.1*cos(118.5)},{4.1*sin(118.5)})};
    \draw [black, thick] plot [smooth] coordinates {(-2.1,2.1) (0,3) (2.1,2.1)};
    \draw [black, thick] plot [smooth] coordinates {(2.1,2.1) (3,1) (-3,1) (-2.1,2.1)};
    \draw[color=black] (-0.2,4.2) node {$x$};
    \draw[color=black] (2.1,1.8) node {$y$};
        \draw[color=black] (-0.15,2.8) node {0};
        \draw[color=black] (-0.15,3.7) node {$\epsilon$};
    \draw[color=black] (-0.2,1.9) node {$-\epsilon$};
    \draw [-to,thick](0,\r1) -- (0,\r2);
    \draw [-to,thick]({\r3*cos(46)},{\r3*sin(46)}) -- ({\r3*cos(45)},{\r3*sin(45)});
    \draw[color=black] (3.8,1.2) node {$\partial\Omega_i$};
    \draw [-to,thick]({2.9*cos(75)},{2.9*sin(75)}) -- ({3.32*cos(74)},{3.32*sin(74)});
    \draw[color=black] (1,2.95) node {$\vec{\boldsymbol{n}}$};
    \end{tikzpicture}
    \caption{An example interface $\partial\Omega_i$ with assigned Fermi coordinates $(x,y)$ where $y\in [0, |\partial\Omega_i| )$ along $\partial\Omega_i$, $x$ is the arclength parameter along the normal lines, with $x=0$ on $\partial\Omega_i$, and $y$ is extended from $\partial \Omega_i$ to be constant along each such line. Here $\vec{\boldsymbol{n}}$ is a normal vector to curve $\partial \Omega_i$.}
    \label{fig:fermi_coords}
\end{figure}

We next use the PDE to obtain regularity for $u$ in the complementary direction. To do this we construct coordinates suitably adapted to each interface and express $T$ in terms of these new coordinates. We consider a unit vector field $\vec{n}$ which is normal to the interface $\partial\Omega_i$. We introduce a coordinate $y\in [0, |\partial\Omega_i|)$ along $\partial\Omega_i$ such that $T(y)=1$. We then introduce Fermi coordinates $(x,y)$ where $x$ is the arclength parameter along the normal lines, with $x=0$ on $\partial\Omega_i$,  and $y$ is extended from $\partial \Omega_i$ to be constant along each such line.  We restrict $x$ to lie in a small interval $(-\epsilon,\epsilon)$ so that no two lines normal to $\partial\Omega_i$ intersect over the intervals $x\in(-\epsilon,\epsilon)$. This can be achieved since $\partial\Omega_i$ is assumed to be $\mathcal{C}^2$. We take $0<\epsilon^\ast\le \epsilon$ such that the lines of constant $y$ value  with $x\in (-\epsilon^\ast,\epsilon^\ast)$ do not intersect any other interface. This can be achieved since all interfaces are mutually disjoint. An example interface with the assigned Fermi coordinates is shown in \Cref{fig:fermi_coords}.

We define the region $\Omega_i^\ast:= \{(x,y)\in\Omega:y\in [0,|\partial\Omega_i|)$, $x\in (-\epsilon^\ast,\epsilon^\ast)\}$.  In these coordinates the Euclidean metric on $\partial\Omega^\ast_i$ 
 is expressed as 
 \begin{equation}
 g_{\mathbbm{E}^2}= dx^2+ g_i(x,y)dy^2,
 \end{equation}
where $g_i(0,y)=1$ and $g_i(x,y)>0$ elsewhere.
 Furthermore, if the interface $\partial\Omega_i\in \mathcal{C}^{k}$, then $g$ is also bounded in $\mathcal{C}^{k-1}(\Omega_i^\ast)$ and its $\mathcal{C}^{k-1}(\Omega_i^\ast)$-norm depends on the $\mathcal{C}^{k}$ norm of $\partial\Omega_i$. 

Then the Laplacian on the region $\Omega_i^\ast$ is 
\begin{equation}\label{eq:Laplacian_fermi}
    \Delta=\partial_{x}^2+\frac{1}{2}\frac{\partial_{x}g_i(x,y)}{g_i(x,y)}\partial_x+\left[g_i(x,y)\right]^{-1}\partial_{y}^2.
\end{equation}

We define $\chi$ to be a smooth cutoff function that is $1$ when $x\in(-\epsilon^\ast/2, \epsilon^\ast/2)$ and $0$ when $x=\pm \epsilon^\ast$. Define the region $\Omega_i^{\ast\ast}:= \{(x,y)\in\Omega:y\in [0,|\partial\Omega_i|)$, $x\in (-\epsilon^\ast/2,\epsilon^\ast/2)\}$. Then we define  $T=\chi\frac{\partial}{\partial y}$.  Invoking the bounds in \cref{eq:bound_nabla_T_j_u}, we have that:

\begin{equation}
    \label{eq:higher norm2}
\partial^j_y(u)\in H^1(\Omega_i^{\ast\ast})
\end{equation}
for all $j\le k$.

Since $\partial^j_y(u)\in H^1(\Omega_i^{\ast\ast})$ we have that $\partial_x\partial^j_y(u)\in L^2(\Omega_i^{\ast\ast})$ for all $j\le k$.

To get bounds on $\partial_x^2\partial^j_y(u)\in L^2(\Omega_i^{\ast\ast})$, we differentiate \cref{eq:Laplacian_fermi} with respect to $y$ a total of $j$ times. We get the equation
\begin{equation}\label{eq:laplace_eqn_isolate_for_dx2u}
    \partial_x^2\partial_y^ju=-\partial_y^j\left[\frac{1}{2}\frac{\partial_{x}g_i(x,y)}{g_i(x,y)}\partial_xu\right]-\partial_y^j\left[\left[g_i(x,y)\right]^{-1}\partial_{y}^2u\right].
\end{equation}
Hence, we can estimate $\partial_x^2\partial_y^ju$ in $L^2(\Omega_i^{\ast\ast})$ by
\begin{equation}\label{eq:dx2dyju_in_L2}
    \|\partial_x^2\partial_y^ju\|_{L^2(\Omega_i^{\ast\ast})}\le\left\lVert\partial_y^j\left[\frac{1}{2}\frac{\partial_{x}g_i(x,y)}{g_i(x,y)}\partial_xu\right]\right\rVert_{L^2(\Omega_i^{\ast\ast})}+\left\lVert\partial_y^j\left[\left[g_i(x,y)\right]^{-1}\partial_{y}^2u\right]\right\rVert_{L^2(\Omega_i^{\ast\ast})}.
\end{equation}

To get higher normal regularity in the region $\Omega_i^{\ast\ast}$, we differentiate \cref{eq:laplace_eqn_isolate_for_dx2u} with respect to $x$ which leads to the equation
\begin{align}\label{eq:laplace_eqn_isolate_for_dx3u}
    \partial^3_x\partial_y^j u&=
    -\frac{1}{2}\partial_y^j\left[\left[g_i(x,y)\right]^{-1}\partial_{x}^2g_i(x,y) u\right]
    +\frac{1}{2}\partial_y^j\left[\left[g_i(x,y)\right]^{-2}\left[\partial_x g_i(x,y)\right]^2 u\right]\nonumber\\
    &\mathrel{\phantom{=}}-\frac{1}{2}\partial_y^j\left[\left[g_i(x,y)\right]^{-1}\partial_xg_i(x,y)\partial_x^2 u\right]
    +\partial_y^j\left[\left[g_i(x,y)\right]^{-2}\partial_xg_i(x,y)\partial_y^{2}u\right]\nonumber\\
    &\mathrel{\phantom{=}}-\partial_y^j\left[\left[g_i(x,y)\right]^{-1}\partial_x\partial_y^{2}u\right]
\end{align}
and the corresponding estimate
\begin{align}
    \left\lVert\partial^3_x\partial_y^j u\right\rVert_{L^2(\Omega_i^{\ast\ast})} &\le
    \left\lVert\frac{1}{2}\partial_y^j\left[\left[g_i(x,y)\right]^{-1}\partial_{x}^2g_i(x,y) u\right]\right\rVert_{L^2(\Omega_i^{\ast\ast})}\nonumber\\
    &\mathrel{\phantom{=}}+\left\lVert\frac{1}{2}\partial_y^j\left[\left[g_i(x,y)\right]^{-2}\left[\partial_x g_i(x,y)\right]^2 u\right]\right\rVert_{L^2(\Omega_i^{\ast\ast})}\nonumber\\
    &\mathrel{\phantom{=}}+\left\lVert\frac{1}{2}\partial_y^j\left[\left[g_i(x,y)\right]^{-1}\partial_xg_i(x,y)\partial_x^2 u\right]\right\rVert_{L^2(\Omega_i^{\ast\ast})}\nonumber\\
    &\mathrel{\phantom{=}}+\left\lVert\partial_y^j\left[\left[g_i(x,y)\right]^{-2}\partial_xg_i(x,y)\partial_y^{2}u\right]\right\rVert_{L^2(\Omega_i^{\ast\ast})}\nonumber\\
    &\mathrel{\phantom{=}}+\left\lVert\partial_y^j\left[\left[g_i(x,y)\right]^{-1}\partial_x\partial_y^{2}u\right]\right\rVert_{L^2(\Omega_i^{\ast\ast})}.
\end{align}

Repeating the same argument as in the concentric circles case in \cref{sm:eq:dr_dtj_u}, we can estimate $T^j\partial_{x} u$ in $H^2(\Omega_i^{\ast\ast})$ by
\begin{align}\label{eq:bound_T_j_dn_u}
    \|T^j\partial_{x}u\|_{H^{2}(\Omega_i^{\ast\ast})}\le\frac{ K_i\cdot\left(4\sup|\nabla T|\cdot \rho_{\sigma_i}\right)^{j+2}}{\min\{\sigma_1,\dots,\sigma_N\}}\|b_1\|_{L^2(\partial\Omega)}
\end{align}
in which $K_i$ is a constant that depends on the trace constant and Poincar\'{e} constant on domain $\Omega_i^{\ast\ast}$ and $\rho_{\sigma_i}$ is $\textrm{max}\{\frac{\sigma_i}{\sigma_{p_i}},\frac{\sigma_{p_i}}{\sigma_{i}}\}$.

\subsection{Regularity of the solution at the interface translates into regularity of the charge densities}

We now derive estimates on the regularity of the charge densities in the case of general domains. \Cref{thm:cgeneral} follows immediately from these estimates.

Using \cref{eq:gammai}, Sobolev embedding, and \cref{eq:bound_T_j_dn_u} we get for $0\le j\le k-2$ and  $0\le\alpha<1$ the bound
\begin{align}
    \|T^j\gamma_i\|_{\mathcal{C}^\alpha(\partial\Omega_i)}&=\left|1-\frac{\sigma_{p_i}}{\sigma_i}\right|\|T^j\partial_{x}u\|_{\mathcal{C}^\alpha(\partial\Omega_i)}\\
    &\le\left|1-\frac{\sigma_{p_i}}{\sigma_i}\right|\|T^j\partial_{x}u\|_{\mathcal{C}^\alpha(\Omega_i^{\ast\ast})}\\
    &\le\left|1-\frac{\sigma_{p_i}}{\sigma_i}\right|c_i\|T^j\partial_{x}u\|_{H^{2}(\Omega_i^{\ast\ast})}\\
    &\le \left|1-\frac{\sigma_{p_i}}{\sigma_i}\right|    \frac{c_i\cdot K_i\cdot \left(4\sup|\nabla T|\cdot\rho_{\sigma_i}\right)^{j+2}}{\min\{\sigma_1,\dots,\sigma_N\}}\|b_1\|_{L^2(\partial\Omega)} \label{eq:Tj_gamma_final_bound}
\end{align}
in which and $\rho_{\sigma_i}$ is $\textrm{max}\{\frac{\sigma_i}{\sigma_{p_i}},\frac{\sigma_{p_i}}{\sigma_{i}}\}$, $c_i$ is the Sobolev embedding constant on $\Omega_i^{\ast\ast}$, and $K_i$ is a constant that depends on the trace constant and Poincar\'{e} constant on domain $\Omega_i^{\ast\ast}$.

\section{Conclusions}
\label{sec:conclusions}

In this paper, we present a novel method for solving the elliptic partial differential equation problem for the electrostatic potential with piecewise constant conductivity. We employ an integral equation approach for which we derive a system of well-conditioned integral equations that can be used to solve the problem. The kernel of the resulting integral operator is smooth.

GMRES is used to solve a matrix vector equation for the charge densities at each grid point. Regarding discretizing the integral equations, we employ two different quadrature schemes. In the case where two interfaces are sufficiently close to one another, we use composite quadrature and split the interface into panels. We then use Gauss-Legendre quadrature on each panel and Lagrange interpolation for interpolating the charge density. For interfaces that are sufficiently far away from all other interfaces we use uniform grid points, trapezoidal rule quadrature, and trigonometric interpolation for interpolating the charge density. The addition of an adaptive method allows our solver to achieve higher accuracy even when the curves are close together at minimal cost to computation time. In the case of interfaces discretized with uniform grid points we use a truncated Fourier series to approximate the scaled charge densities. On the other hand, in the case of interfaces discretized with composite quadrature (panels) we use a truncated Legendre series polynomial to approximate the scaled charge densities. In either case, if either of the two highest mode coefficients are above a specified threshold, then the interface is refined.

When the evaluation point is not near or on an interface or the boundary, we evaluate the single layer potentials using the same quadrature as used to solve the integral equations. When the evaluation point is near or on an interface or the boundary, we approximate interfaces or the boundary with line segments near the evaluation point and employ an analytic expression for the single layer potential due to a line segment.

To illustrate the effectiveness of our method, we solve the elliptic partial differential equation problem for four test cases. Our method is compared against a popular open-source platform and is shown to be superior since it produces more accurate results in less computational time. In addition, our method is also shown to easily handle problems of increasingly higher complexity involving up to 155 different regions of constant conductivity.

Numerous avenues of future work exist. One could adapt our method to solve the problem in which the potential $u$ has a prescribed jump across the interfaces by including double layer potentials. Also, the method could be adapted to the three dimensional version of the problem. Each of these generalizations would enable several applications of the method.

\bibliographystyle{siamplain}
\bibliography{references}

\makeatletter\@input{supplementx.tex}\makeatother

\end{document}


\maketitle

\section{Existence and Regularity}
\label{sm:sec:Existence and Regularity}

In this section, we will provide a theoretical study of the system of integral equations \crefrange{eq:inteqn1}{eq:inteqn2} in the special case where the regions of constant conductivity are defined by nested concentric circles.

\begin{theorem}
\label{sm:thm:concentric}
Consider the solution of problem \crefrange{eq:inteqn1}{eq:inteqn2}, with $\partial\Omega_i$ \\nested concentric circles inside the unit disk $\mathbbm{D}$ and injected current $b_1\in H^k(\partial\mathbbm{D})$. Consider the charge densities $\gamma_i(\cdot)$ at each of the inner interfaces. Then these functions are ${\cal C}^{k-2,\alpha}$ regular for $0\le\alpha<1$, and their norms are bounded as follows: ($\theta$ is the standard angular polar coordinate)
\begin{equation}
\|\gamma_i(\theta)\|_{{\cal C}^{k-2,\alpha}(\partial\Omega_i)}\le
\left|1-\frac{\sigma_{p_i}}{\sigma_i}\right|\frac{M}{\min\{\sigma_1,\dots,\sigma_N\}}\frac{6}{(r_N)^2}\|b_1\|_{H^k(\partial\mathbbm{D})}
\end{equation}
in which $r_N$ is the radius of the innermost circle and $M$ is a constant that depends on the constants in the Sobolev embedding, the trace theorem, and the Poincar\'{e} inequality on domain $\mathbbm{D}$.
\end{theorem}

The constant $M$ does not blow up as two interfaces (circles in this case) approach each other; of course on the other hand, the bounds are in terms of a stronger norm of the injected current than in the general case of \cref{thm:cgeneral}.

\subsection{Regularity of the solution at the interfaces, in the context of concentric circles}

\begin{figure}[ht]
\centering
\begin{tikzpicture}[scale=1, transform shape][line cap=round,line join=round,>=triangle 45,x=1cm,y=1cm]
\draw [line width=1pt,color=black] (0,0) circle[radius=2cm];
\draw [line width=1pt,color=black] (0,0) circle[radius=1.4cm];
\draw [line width=1pt,color=black] (0,0) circle[radius=1cm];
\draw [line width=1pt,color=black] (0,0) circle[radius=0.6cm];
\draw [line width=0.75pt,color=black,fill=black] (-0.65,0.65) circle[radius=0.01cm];
\draw [line width=0.75pt,color=black,fill=black] (-0.575,0.575) circle[radius=0.01cm];
\draw [line width=0.75pt,color=black,fill=black] (-0.5,0.5) circle[radius=0.01cm];
\draw [line width=0.5pt,color=black] (0,0) -- (0.5638,0.2052);
\draw [line width=0.5pt,color=black] (0,0) -- (1.0724,0.8999);
\draw [line width=0.5pt,color=black] (0,0) -- (1,1.732);
\begin{scriptsize}
\draw[color=black] (-1.2,1.2) node {$\sigma_1$};
\draw[color=black] (-0.83,0.83) node {$\sigma_2$};
\draw[color=black] (-0.25,0.25) node {$\sigma_N$};
\draw[color=black] (0.8,0.15) node {$r_N$};
\draw[color=black] (1.3,0.9) node {$r_2$};
\draw[color=black] (1.2,1.8) node {$r_1$};
\end{scriptsize}
\end{tikzpicture}
\caption{Domain for the problem with $N$ concentric circular regions of constant conductivity. The conductivity jumps across the outer boundary of region $i\ne1$, a circle of radius $r_i$, from $\sigma_i$ to $\sigma_{i-1}$.} \label{sm:fig:Domain2ConcentircConductivities}
\end{figure}
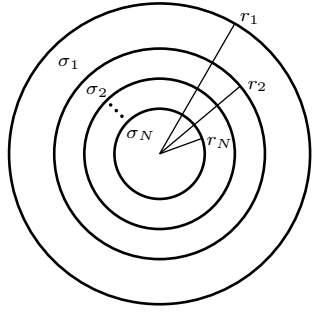

We consider the special case where the regions of piecewise constant conductivity are defined by $N\in\mathbbm{N}$ concentric circles, of radii $r_1=1>r_2>\dots>r_N>0$. \Cref{sm:fig:Domain2ConcentircConductivities} shows the domain and conductivities that we consider. 

We denote the unit disk by $\mathbbm{D}$ and use the standard polar coordinates $(r,\theta)$. Let $b_1\in H^k(\partial\mathbbm{D})$ for some $k\in\mathbbm{N}$. The weak formulation of \cref{eq: PDE,eq: BC} is
\begin{equation}
\label{sm:eq: weak formulation}
    \int_\mathbbm{D} \sigma \langle \nabla u, \nabla \phi \rangle~dA = \int_{\partial \mathbbm{D}} b_1 \phi~dl
\end{equation}
for all $\phi\in \mathcal{C}^\infty$ on the closed disk.

Since $\sigma\in L^\infty$, from general PDE theory we obtain that there exists a unique weak solution $u\in H^{1}_{0}(\mathbbm{D})$~\cite{evans10dup}. Setting $\phi=u$ and using the Cauchy–Schwarz inequality, the trace theorem, and the Poincar\'{e} inequality we obtain
\begin{align}
    \min\{\sigma_1,\dots,\sigma_N\}\int_\mathbbm{D}|\nabla u|^2 dA &\le \int_\mathbbm{D}\sigma |\nabla u|^2 dA \\
    &=\int_{\partial \mathbbm{D}} b_1 u~dl\\
    &\le \frac{1}{2}\mu\int_{\partial \mathbbm{D}}b_1^2~dl+\frac{1}{2}\mu^{-1}\int_{\partial \mathbbm{D}}u^2~dl\\
    &\le \frac{1}{2}\mu\int_{\partial \mathbbm{D}}b_1^2~dl+\frac{1}{2}\mu^{-1}c_1^2\|u\|^2_{H^{1}(\mathbbm{D})}\\
    &= \frac{1}{2}\mu\int_{\partial \mathbbm{D}}b_1^2~dl+\frac{1}{2}\mu^{-1}c_1^2\left(\int_{\mathbbm{D}}|u|^2+|\nabla u|^2~dA\right)\\
    &\le \frac{1}{2}\mu\int_{\partial \mathbbm{D}}b_1^2~dl+\frac{1}{2}\mu^{-1}c^2\int_{\mathbbm{D}}|\nabla u|^2~dA
\end{align}
in which $c_1$ is the trace constant and $\sqrt{\frac{c^2}{c_1^2}-1}$ is the Poincar\'{e} constant for the domain $\mathbbm{D}$.

Thus, choosing $\mu=\frac{c^2}{\min\{\sigma_1,\dots,\sigma_N\}}$ we get the bound
\begin{align}\label{sm:eq:bound_on_grad_u}
    \|\nabla u\|_{L^2(\mathbbm{D})} &\le \frac{c}{\min\{\sigma_1,\dots,\sigma_N\}}\|b_1\|_{L^2(\partial\mathbbm{D})}.
\end{align}
in which $c$ is a constant depending on the trace constant and Poincar\'{e} constant for the domain $\mathbbm{D}$.

We seek to prove better regularity for $u$. First, we examine tangential regularity. Let $w\in \mathcal{C}^\infty$ and take $\phi=\partial_\theta^j w$ for some $j\le k$ in \cref{sm:eq: weak formulation} \textit{i}.\textit{e}.,
\begin{equation}
        \int_\mathbbm{D} \sigma \langle \nabla u, \nabla \partial_\theta^j w \rangle~dA = \int_{\partial \mathbbm{D}} b_1 \partial_\theta^j w~dl.
\end{equation}
We integrate by parts in $\theta$ to give
\begin{equation}
        \int_\mathbbm{D} \sigma \langle \nabla \partial_\theta^j u, \nabla w \rangle~dA = \int_{\partial \mathbbm{D}} \partial_\theta^j b_1  w~dl
\end{equation}
so $\partial_\theta^j u$ is a weak solution to the same PDE problem \cref{eq: PDE} with Neumann boundary condition $\sigma\partial_{\boldsymbol{n}}u = \partial_\theta^j b_1$. Furthermore, $\partial_\theta^j u\in H^1(\mathbbm{D})$ for $j\le k$ with the estimate (note that $\int_{\mathbbm{D}}\partial_\theta u~dA=0$)
\begin{equation}\label{sm:eq:bound_on_grad_partial_jt_u}
    \|\nabla \partial_\theta^j u\|_{L^2(\mathbbm{D})} \le \frac{c}{\min\{\sigma_1,\dots,\sigma_N\}}\|\partial_\theta^j b_1\|_{L^2(\partial\mathbbm{D})}
\end{equation}
in which $c$ is a constant depending on the trace constant and Poincar\'{e} constant for the domain $\mathbbm{D}$.

Next, we examine normal regularity in the annulus $\Omega_{n}\setminus\Omega_{n+1}$ up to the interface for $n\in\{1,\dots,N-1\}$. Since $\partial_\theta^j u\in H^1(\mathbbm{D})$ we get that $\partial_r\partial_\theta^j u\in L^2(\Omega_{n}\setminus\Omega_{n+1})$ and $\partial_\theta^{j+1} u\in L^2(\Omega_{n}\setminus\Omega_{n+1})$ for $j\le k$. We use the PDE \cref{eq: PDE} to get bounds on $\partial_r^i\partial_\theta^j u$ for $i+j\le k+1$ and $i\le3$. In the annulus $\Omega_{n}\setminus\Omega_{n+1}$, $\partial_\theta^j u$ is harmonic so
\begin{equation}\label{sm:eq:laplace_eqn_isolate_for_dr2u}
    \partial^2_r\partial_\theta^j u=-\frac{1}{r}\partial_r\partial_\theta^j u-\frac{1}{r^2}\partial_\theta^{j+2} u.
\end{equation}
Hence, we can estimate $\partial_r^2 \partial_\theta^j u$ in $L^2(\Omega_{n}\setminus\Omega_{n+1})$ by
\begin{align}
    \left\lVert\partial^2_r \partial_\theta^j u\right\rVert_{L^2(\Omega_{n}\setminus\Omega_{n+1})}&\le\frac{1}{r_{n+1}}\left\lVert\partial_r \partial_\theta^j u\right\rVert_{L^2(\Omega_{n}\setminus\Omega_{n+1})}+\left(\frac{1}{r_{n+1}}\right)^2\left\lVert \partial_\theta^{j+2} u\right\rVert_{L^2(\Omega_{n}\setminus\Omega_{n+1})}\\
    &\le\frac{1}{r_{n+1}}\left\lVert\nabla \partial_\theta^j u\right\rVert_{L^2(\Omega_{n}\setminus\Omega_{n+1})}+\left(\frac{1}{r_{n+1}}\right)^2\left\lVert \nabla\partial_\theta^{j+1} u\right\rVert_{L^2(\Omega_{n}\setminus\Omega_{n+1})}
\end{align}
for $j\le k-1$. To get higher normal regularity in the annulus $\Omega_{n}\setminus\Omega_{n+1}$, we differentiate \cref{sm:eq:laplace_eqn_isolate_for_dr2u} with respect to $r$ which leads to the equation
\begin{equation}\label{sm:eq:laplace_eqn_isolate_for_dr3u}
    \partial^3_r\partial_\theta^j u=\frac{1}{r^2}\partial_r\partial_\theta^j u-\frac{1}{r}\partial_r^2\partial_\theta^j u+\frac{2}{r^3}\partial_\theta^{j+2}u-\frac{1}{r^2}\partial_r\partial_\theta^{j+2} u
\end{equation}
and the corresponding estimate
\begin{align}
    \left\lVert\partial^3_r\partial_\theta^j u\right\rVert_{L^2(\Omega_{n}\setminus\Omega_{n+1})}&\le\left(\frac{1}{r_{n+1}}\right)^2\left\lVert\partial_r\partial_\theta^j u\right\rVert_{L^2(\Omega_{n}\setminus\Omega_{n+1})}+\frac{1}{r_{n+1}}\left\lVert\partial_r^2\partial_\theta^j u\right\rVert_{L^2(\Omega_{n}\setminus\Omega_{n+1})}\nonumber\\
    &\mathrel{\phantom{=}}+2\left(\frac{1}{r_{n+1}}\right)^3\left\lVert\partial_\theta^{j+2}u\right\rVert_{L^2(\Omega_{n}\setminus\Omega_{n+1})}\nonumber\\
    &\mathrel{\phantom{=}}+\left(\frac{1}{r_{n+1}}\right)^2\left\lVert\partial_r\partial_\theta^{j+2} u\right\rVert_{L^2(\Omega_{n}\setminus\Omega_{n+1})}\\
    &\le 2\left(\frac{1}{r_{n+1}}\right)^2\left\lVert\nabla \partial_\theta^j u\right\rVert_{L^2(\Omega_{n}\setminus\Omega_{n+1})} \\
    &\mathrel{\phantom{=}}+3\left(\frac{1}{r_{n+1}}\right)^2\left\lVert \nabla\partial_\theta^{j+1} u\right\rVert_{L^2(\Omega_{n}\setminus\Omega_{n+1})} \\
    &\mathrel{\phantom{=}}+\left(\frac{1}{r_{n+1}}\right)^2\left\lVert \nabla\partial_\theta^{j+2} u\right\rVert_{L^2(\Omega_{n}\setminus\Omega_{n+1})}
\end{align}
for $j\le k-2$. 

We can estimate $\partial_r\partial_\theta^j u$ in $H^2$ on the annulus $\Omega_{n}\setminus\Omega_{n+1}$ by
\begin{align}
    \left\lVert \partial_r\partial_\theta^j u \right\rVert_{H^2(\Omega_{n}\setminus\Omega_{n+1})}&= \left\lVert \partial_r\partial_\theta^j u \right\rVert_{L^2(\Omega_{n}\setminus\Omega_{n+1})} + \left\lVert \partial_r\partial_\theta^{j+1} u \right\rVert_{L^2(\Omega_{n}\setminus\Omega_{n+1})} \nonumber\\
    &\mathrel{\phantom{=}}+\left\lVert \partial_r\partial_\theta^{j+2} u \right\rVert_{L^2(\Omega_{n}\setminus\Omega_{n+1})}+\left\lVert \partial^2_r\partial_\theta^j u \right\rVert_{L^2(\Omega_{n}\setminus\Omega_{n+1})} \nonumber\\
    &\mathrel{\phantom{=}}+\left\lVert \partial^2_r\partial_\theta^{j+1} u \right\rVert_{L^2(\Omega_{n}\setminus\Omega_{n+1})}+\left\lVert \partial^3_r\partial_\theta^j u \right\rVert_{L^2(\Omega_{n}\setminus\Omega_{n+1})} \\
    &\le \left\lVert \nabla\partial_\theta^j u \right\rVert_{L^2(\Omega_{n}\setminus\Omega_{n+1})} + \left\lVert \nabla\partial_\theta^{j+1} u \right\rVert_{L^2(\Omega_{n}\setminus\Omega_{n+1})} \nonumber\\
    &\mathrel{\phantom{=}}+\left\lVert \nabla\partial_\theta^{j+2} u \right\rVert_{L^2(\Omega_{n}\setminus\Omega_{n+1})}+\left\lVert \partial^2_r\partial_\theta^j u \right\rVert_{L^2(\Omega_{n}\setminus\Omega_{n+1})} \nonumber\\
    &\mathrel{\phantom{=}}+\left\lVert \partial^2_r\partial_\theta^{j+1} u \right\rVert_{L^2(\Omega_{n}\setminus\Omega_{n+1})}+\left\lVert \partial^3_r\partial_\theta^j u \right\rVert_{L^2(\Omega_{n}\setminus\Omega_{n+1})} \\
    &\le \left(1+\frac{1}{r_{n+1}}+2\left(\frac{1}{r_{n+1}}\right)^2\right)\left\lVert \nabla\partial_\theta^j u \right\rVert_{L^2(\Omega_{n}\setminus\Omega_{n+1})} \nonumber\\
    &\mathrel{\phantom{=}}+\left(1+\frac{1}{r_{n+1}}+4\left(\frac{1}{r_{n+1}}\right)^2\right) \left\lVert \nabla\partial_\theta^{j+1} u \right\rVert_{L^2(\Omega_{n}\setminus\Omega_{n+1})} \nonumber\\
    &\mathrel{\phantom{=}}+\left(1+2\left(\frac{1}{r_{n+1}}\right)^2\right)\left\lVert \nabla\partial_\theta^{j+2} u \right\rVert_{L^2(\Omega_{n}\setminus\Omega_{n+1})}
\end{align}

We can estimate $\partial_r\partial_\theta^j u$ in $H^2$ on the annulus $\Omega\setminus\Omega_{N}$ by
\begin{align}
    \left\lVert \partial_r\partial_\theta^j u \right\rVert_{H^2(\Omega\setminus\Omega_{N})}&\le\frac{6}{(r_N)^2}\left(\left\lVert \nabla\partial_\theta^{j} u \right\rVert_{L^2(\Omega\setminus\Omega_{N})}+\left\lVert \nabla\partial_\theta^{j+1} u \right\rVert_{L^2(\Omega\setminus\Omega_{N})}\right.\nonumber\\
    &\mathrel{\phantom{=}}\qquad\qquad\left.+\left\lVert \nabla\partial_\theta^{j+2} u \right\rVert_{L^2(\Omega\setminus\Omega_{N})}\right)\label{sm:eq:dr_dtj_u}
\end{align}

\subsection{Regularity of the solution at the interface translates into regularity of the charge densities}

We now derive estimates on the regularity of the charge densities in the case of concentric circle domains. \cref{sm:thm:concentric} follow immediately from these estimates.

Using \cref{eq:gammai}, Sobolev embedding, \cref{sm:eq:bound_on_grad_partial_jt_u}, and \cref{sm:eq:dr_dtj_u}, we get for $0\le j\le k-2$ and $0\le\alpha<1$ the bound
\begin{align}
    \|\partial_\theta^j\gamma_i\|_{{\cal C}^\alpha(\partial\Omega_i)}&=\left|1-\frac{\sigma_{p_i}}{\sigma_i}\right|\|\partial_{r}\partial_\theta^ju\|_{{\cal C}^\alpha(\partial\Omega_i)}\\
    &\le\left|1-\frac{\sigma_{p_i}}{\sigma_i}\right|\|\partial_{r}\partial_\theta^ju\|_{{\cal C}^\alpha\left(\mathbbm{D}\setminus\left(D\left(\boldsymbol{0},\frac{1}{4}\right)\cap\Omega_N\right)\right)}\\
    &\le \left|1-\frac{\sigma_{p_i}}{\sigma_i}\right|c\|\partial_r \partial_\theta^ju\|_{H^{2}\left(\mathbbm{D}\setminus\left(D\left(\boldsymbol{0},\frac{1}{4}\right)\cap\Omega_N\right)\right)}\\
    &\le\left|1-\frac{\sigma_{p_i}}{\sigma_i}\right| \frac{c\cdot K}{\min\{\sigma_1,\dots,\sigma_N\}}\frac{6}{(r_N)^2}\| b_1\|_{H^{j+2}(\partial\mathbbm{D})}
\end{align}
in which $c$ is the Sobolev embedding constant on on disks of size $\frac{1}{2}$ (which is uniformally bounded), $K$ is a constant that depends on the trace constant and Poincar\'{e} constant on domain $\mathbbm{D}$, and $r_N$ is the radius of the innermost circle. Here $D\left(\boldsymbol{0},\frac{1}{4}\right)$ is the disk of radius $\frac{1}{4}$ centered at the origin.

Remark:  We note that in the special case of concentric circles and for injected currents that are smooth enough as functions of $\theta$ on the outer boundary, each extra derivative of our solution is bounded without the bounds deteriorating if the two circles are close together. In contrast, in the general case, if two interfaces are close, then the term $\sup|\nabla T|$ in \cref{eq:Tj_gamma_final_bound}, for the domain we constructed will necessarily become large. Moreover, an interface with large geodesic curvature $\kappa$ will also cause the constants in the higher-derivatives estimates to deteriorate (\textit{i}.\textit{e}., become large). Thus, the more general estimate is more versatile, yet more sensitive to the interfaces being either mutually close, or very curved. 

\bibliographystyle{siamplain}
\bibliography{references}

\makeatletter\@input{articlex.tex}\makeatother